\documentclass[10pt]{svjour3}                     
\smartqed  
\usepackage{amsmath,amssymb} 
\usepackage{url,hyperref}
\usepackage{mathtools}
\usepackage{multirow}
\usepackage{graphicx}
\usepackage{color}
\usepackage{algorithm}

\makeatletter
\renewcommand\@biblabel[1]{#1.}
\makeatother

\newcommand{\Id}{\mathbb{I}}
\newcommand{\R}{\mathbb{R}}

\newcommand{\set}[1]{\left\{#1\right\}}
\newcommand{\sets}[1]{\{#1\}}

\newcommand{\norms}[1]{\Vert#1\Vert}

\newcommand{\iprods}[1]{\langle #1\rangle}
\newcommand{\Eproof}{\hfill $\square$}

\newcommand{\dom}[1]{\mathrm{dom}(#1)}
\newcommand{\zero}[1]{{\boldsymbol{0}}}

\newcommand{\zer}[1]{\mathrm{zer}(#1)}
\newcommand{\gra}[1]{\mathrm{gra}(#1)}

\newcommand{\Hc}{\mathcal{H}}

\newcommand{\Lc}{\mathcal{L}}
\newcommand{\Qc}{\mathcal{Q}}
\newcommand{\Tc}{\mathcal{T}}

\newcommand{\Vc}{\mathcal{V}}

\newcommand{\Ec}{\mathcal{E}}
\newcommand{\Nc}{\mathcal{N}}

\newcommand{\BigO}[1]{\mathcal{O}\left(#1\right)}
\newcommand{\SmallO}[1]{o\left(#1\right)}

\newcommand{\beforesubsec}{\vspace{-3.75ex}}
\newcommand{\aftersubsec}{\vspace{-2.2ex}}
\newcommand{\beforesec}{\vspace{-3.5ex}}
\newcommand{\aftersec}{\vspace{-2.25ex}}

\newcommand{\revise}[1]{\textcolor{black}{#1}}

\newcommand{\myeq}[2]{\vspace{-0.5ex}\begin{equation}\label{#1}#2\vspace{-0.5ex}\end{equation}}
\newcommand{\myeqn}[1]{\vspace{-0.25ex}\begin{equation*}#1\vspace{-0.5ex}\end{equation*}}

\journalname{Comput Optim Appl}

\begin{document}

\title{From Halpern's Fixed-Point Iterations to Nesterov's Accelerated Interpretations for Root-Finding Problems
}

\titlerunning{From Halpern's Fixed-Point Iterations to Nesterov's Accelerated Interpretations}        

\author{Quoc Tran-Dinh
}


\institute{Quoc Tran-Dinh \at
		Department of Statistics and Operations Research,  The University of North Carolina at Chapel Hill, 333 Hanes Hall, CB\# 3260, UNC Chapel Hill, NC 27599-3260, USA\\ 
		Email: \url{quoctd@email.unc.edu}
}

\date{Received: date / Accepted: date}

\maketitle

\begin{abstract}
We derive an equivalent form of Halpern's fixed-point iteration scheme for solving a co-coercive equation (also called a root-finding problem), which can be viewed as a Nesterov's accelerated interpretation.
We show that one method is equivalent to another via a simple transformation, leading to a straightforward  convergence proof for Nesterov's accelerated scheme.
Alternatively, we directly establish convergence rates of Nesterov's accelerated variant, and as a consequence, we obtain a new convergence rate of Halpern's fixed-point iteration.
Next, we apply our results to different methods to solve monotone inclusions, where our convergence guarantees are applied.
Since the gradient/forward scheme requires the co-coerciveness of the underlying operator, we derive new Nesterov's accelerated variants for both recent extra-anchored gradient  and past-extra anchored gradient methods in the literature. 
These variants alleviate  the co-coerciveness condition by only assuming the monotonicity and Lipschitz continuity of the underlying operator.
Interestingly, our new Nesterov's accelerated interpretation of the past-extra anchored gradient method involves two past-iterate correction terms.
This formulation is expected to guide us developing new Nesterov's accelerated methods for minimax problems and their continuous views without co-coericiveness. 
We test our theoretical results on two numerical examples, where the actual convergence rates match well the theoretical ones up to a constant factor.
\end{abstract}


\beforesec
\section{Introduction}\label{sec:intro}
\aftersec
Approximating a solution of a maximally monotone inclusion is a fundamental problem in optimization, nonlinear analysis, mechanics, and machine learning, among many other areas, see, e.g.,  \cite{Bauschke2011,reginaset2008,Facchinei2003,phelps2009convex,Rockafellar2004,Rockafellar1976b,ryu2016primer}.
This problem lies at the heart of monotone operator theory, and has been intensively studied in the literature for many decades, see, e.g., \cite{Bauschke2011,Lions1979,Rockafellar1976b,ryu2016primer} as a few examples.
Various numerical methods, including proximal-point-type, gradient/forward, extragradient, past-extragradient, and their variants  have been proposed to solve this problem, and its extensions as well as special cases \cite{Bauschke2011,Davis2015,Facchinei2003,Monteiro2010a,Nemirovskii2004,popov1980modification,Rockafellar1976b}.
When the underlying operator is the sum of two or multiple maximally monotone operators, forward-backward splitting, forward-backward-forward splitting, Douglas-Rachford splitting, projective splitting methods, and their variants have been extensively developed for approximating solutions of this problem under different assumptions and context, see, e.g., \cite{Bauschke2011,Combettes2005,Davis2015,Lions1979,malitsky2015projected,popov1980modification,tseng2000modified} as a few references.

\textbf{Motivation and related work.}
In the last decades, accelerated first-order methods have become an extremely popular and attractive research topic in optimization and related fields due to their applications to large-scale optimization problems in machine learning, statistics, signal and image processing, and engineering, see, e.g., \cite{Beck2009,Nesterov1983,Nesterov2004,Nesterov2005c}. 
In this research theme, Nesterov's accelerated approach \cite{Nesterov1983} presents as a leading research topic for many years, and remains emerging in optimization community.
This well-known technique has been extended to different directions, including minimax problems, variational inequalities (VIPs), and monotone inclusions \cite{attouch2020convergence,bot2022bfast,kim2021accelerated,mainge2021fast}.
Convergence rates  of these methods have been intensively studied, which show significant improvements from $\BigO{1/k}$ to $\BigO{1/k^2}$ rates, where $k$ is the iteration counter.
The latter rate matches the lower bound rates in different settings using different criteria, see, e.g., \cite{kim2016optimized,Nesterov2004,ouyang2018lower}.
In recent years, many papers, including \cite{attouch2020convergence,kim2021accelerated,mainge2021fast}, have focused on developing Nesterov's accelerated schemes for monotone inclusions.
They have proven $\BigO{1/k^2}$-convergence rates, and also $\SmallO{1/k^2}$ convergence rates on the squared norm of the residual operator associated with the problem, while obtaining asymptotic convergence on iterate sequences, see \cite{attouch2016rate,chambolle2015convergence,kim2016optimized}. 
Note that the problem of approximating a solution of a maximally monotone inclusion can be reformulated equivalently to approximating a fixed-point of a nonexpansive operator \cite{Bauschke2011}. 
Therefore, theory and solution methods from one field can be applied to the other and vice versa. 

Alternatively, Halpern's fixed-point iteration is a classical method in fixed-point theory rooted from \cite{halpern1967fixed} to approximate a fixed-point of a nonexpansive operator, see  \cite{bauschke1996approximation,kornlein2015quantitative,wittmann1992approximation}.
This method has recently attracted great attention due to its ability to accelerate convergence rate in terms of operator residual norm.
Lieder specifically proved an $\BigO{1/k^2}$ rate on the squared norm of the operator residual for Halpern's fixed-point iteration in  \cite{lieder2021convergence}, but \cite{sabach2017first} appears to be the first one achieving this rate for a variant of Halpern's fixed-point method.
Unlike Nesterov's accelerated method which is originally developed for solving convex optimization problems and its convergence rate is given in terms of the objective residual in the general convex case, Halpern's fixed-point iteration is  proposed to approximate a fixed-point of a nonexpansive operator, which is much more general than convex optimization, and hence, very convenient to extend to maximally monotone inclusions, and, in particular, minimax problems, game theory, robust optimization, online learning, and reinforcement learning, see, e.g., \cite{diakonikolas2020halpern,lee2021fast,yoon2021accelerated}. 

A natural question is therefore arising: \textit{What is the relation between the two accelerated techniques?}
Such a type of questions was mentioned in \cite{yoon2021accelerated}.
Since both schemes come from different roots, at first glance, it is unclear to see a close relation between Nesterov's accelerated and the Halpern fixed-point schemes.
Nesterov's accelerated method is perhaps rooted from the gradient descent scheme in convex optimization with an additional momentum term  \cite{Polyak1964}.
\revise{Its acceleration behavior has been explained through different view points, including geometric interpretation \cite{bubeck2015geometric} and continuous views via ordinary differential equations (ODEs) and variational perspectives \cite{attouch2022ravine,attouch2019convergence,shi2021understanding,Su2014,wibisono2016variational}.
In addition, Nesterov's accelerated method has various variants \cite{attouch2019convergence,kim2021accelerated,kim2016optimized,mainge2021fast,shi2021understanding}. 
For instance, \cite{kim2016optimized} derived an ``optimized'' Nesterov's accelerated variant to solve composite convex optimization problems which is slightly different from the original one in \cite{Nesterov2004}.
As other examples, both \cite{kim2021accelerated} and \cite{mainge2021fast} proposed  Nesterov's accelerated schemes using different ``correction'' terms to solve monotone inclusions.
}
Alternatively, Halpern's fixed-point method was proposed in \cite{halpern1967fixed} since 1967, existing convergence guarantees are essentially asymptotic or slow convergence rates \cite{wittmann1992approximation}.
Its accelerated rate has just recently been established in \cite{lieder2021convergence,sabach2017first} and followed up by a number of works, including \cite{diakonikolas2020halpern,lee2021fast,yoon2021accelerated}.
Interestingly, the analysis of both schemes is quite different, but still relies on appropriate Lyapunov or potential functions, and variable parameters (e.g., stepsizes, extrapolation, and momentum parameters).

\vspace{0.5ex}
\textbf{Our contribution.}
In this paper, we show that Halpern's fixed-point method can be transformed into Nesterov's accelerated interpretation and vice versa. 
We first present our results on approximating a solution of a co-coercive equation, and then extend them to other schemes.
In the first case, we establish that the iterate sequences generated by both schemes are identical, but the choice of parameters in these schemes could be different, leading to different convergence guarantees.
In fact, we can obtain the convergence guarantee  of one scheme from another and vice versa. 
Then, by utilizing our analysis, we prove that a number of methods, including proximal-point, forward-backward splitting, Douglas-Rachford splitting, and three-operator splitting schemes can be easily accelerated and achieve faster rates compared to their classical counterparts.

Note that, in convex optimization, we often use the objective residual as a potential energy term to form a Lyapunov or an energy function for establishing convergence rates. 
Moreover, this objective residual presents as a main metric to measure the approximate optimality of the current iterate.
However, such an objective function does not exist in root-finding problems. 
Therefore, establishing $\BigO{1/k^2}$ and $\SmallO{1/k^2}$ convergence rates for a root-finding problem requires different metric as well as different techniques  than those in convex optimization. 

Our next contribution is to show that the recent extra-anchored gradient (EAG) proposed in \cite{lee2021fast,yoon2021accelerated} can be transformed into Nesterov's accelerated interpretation.
We provide an alternative analysis using a slightly different Lyapunov function and still obtain the same convergence rates as in \cite{lee2021fast,yoon2021accelerated}.
One important fact is that this accelerated method works with a monotone and Lipschitz continuous operator instead of a co-coercive one as in \cite{attouch2020convergence,kim2021accelerated,mainge2021fast}.
This approach is expected to provide an initial step toward understanding Nesterov's accelerated behavior on non-co-coercive operators, and possibly their continuous views.

Finally, we derive Nesterov's accelerated variant of the past-extra anchored gradient method in \cite{tran2021halpern} and provide a different convergence rate analysis than that of  \cite{tran2021halpern} by using two different stepsizes for the extra-gradient step.
Interestingly, such a new scheme is very different from existing ones, e.g., \cite{attouch2020convergence,kim2021accelerated,mainge2021fast}, due to the use of two consecutive past iterates in the momentum/correction terms.
To the best of our knowledge, this algorithm is the first one illustrating that we can build Nesterov's accelerated scheme for monotone equations without co-coerciveness. 

While we are working on this paper, a recent manuscript \cite{park2022exact} comes to our attention, which closely relates to Subsection~\ref{subsec:NAG2HP} of this paper.
However, \cite{park2022exact} mainly focuses on establishing exact optimal rates of accelerated schemes for maximally monotone inclusions, whereas we study the relation between the two accelerated approaches and their convergence guarantees for   different methods, including extra-anchored gradient.
In addition, after our manuscript was uploaded online on \href{researchgate.net}{researchgate.net}, another work \cite{bot2022fast}, uploaded a few months later on arxiv, also studies Nesterov's accelerated variants for both extragradient and past-extragradient methods, but from a discretization of  dynamical systems.

\vspace{0.5ex}
\textbf{Paper organization.}
The rest of this paper is organized as follows.
Section~\ref{sec:background} provides some necessary background of monotone operators, which will be used in the sequel.
Section~\ref{sec:main_part1} presents the equivalence between Halpern's fixed-point iteration and Nesterov's accelerated scheme for solving a co-coercive equation.
Next, we discuss its application to other methods
 in Section~\ref{sec:mono_inclusion}.
Then, Section~\ref{sec:EAG_sec} investigates the connection between extra-anchored gradient method and its variants and the corresponding Nesterov's interpretations.
Finally, Section~\ref{sec:num_experiments} provides two numerical examples to verify our theoretical results.

\beforesec
\section{Background of Monotone Operators}\label{sec:background}
\aftersec
%
We work with finite-dimensional spaces $\R^p$ and $\R^n$ endowed  with the standard inner product $\iprods{\cdot, \cdot}$ and Euclidean norm $\norms{\cdot}$.
For a set-valued mapping $G : \R^p \rightrightarrows 2^{\R^p}$, $\dom{G} = \{ x \in\R^p : Gx \not= \emptyset \}$ denotes its domain, $\gra{G} = \{ (x, y) \in \R^p\times \R^p : y \in Gx \}$ denotes its graph, where $2^{\R^p}$ is the set of all subsets of $\R^p$.
The inverse of $G$ is defined as $G^{-1}y := \sets{x \in \R^p : y \in Gx}$.

\textbf{Monotonicity.}
For a set-valued mapping $G : \R^p \rightrightarrows 2^{\R^p}$, we say that $G$ is monotone if $\iprods{u - v, x - y} \geq 0$ for all $(x, u), (y, v) \in \gra{G}$.
$G$ is said to be $\mu_G$-strongly monotone (also called coercive) if $\iprods{u - v, x - y} \geq \mu_G\norms{x - y}^2$ for all $(x, u), (y, v) \in \gra{G}$, where $\mu_G > 0$ is called the strong monotonicity parameter.
If $G$ is single-valued, then these conditions reduce to $\iprods{Gx - Gy, x - y} \geq 0$ and $\iprods{Gx - Gy, x - y} \geq \mu_G\norms{x - y}^2$ for all $x, y\in\dom{G}$, respectively.
We say that $G$ is maximally monotone if $\gra{G}$ is not properly contained in the graph of any other monotone operator.
Further details can be found, e.g., in \cite{Bauschke2011,ryu2016primer}.

\textbf{Lipschitz continuity and co-coerciveness.}
A single-valued operator $G$ is said to be $L$-Lipschitz continuous if $\norms{Gx - Gy} \leq L\norms{x - y}$ for all $x, y\in\dom{G}$, where $L \geq 0$ is the Lipschitz constant. 
If $L = 1$, then we say that $G$ is nonexpansive, while if $L \in [0, 1)$, then we say that $G$ is $L$-contractive, and $L$ is its contraction factor.
We say that $G$ is $\frac{1}{L}$-co-coercive if $\iprods{Gx - Gy, x - y} \geq \frac{1}{L}\norms{Gx - Gy}^2$ for all $x, y\in\dom{G}$.
If $L = 1$, then we say that $G$ is firmly nonexpansive.
Note that if $G$ is $\frac{1}{L}$-co-coercive, then it is also monotone and $L$-Lipschitz continuous, but the reverse is not true in general.
If $L < 0$, then we say that $G$ is $\frac{1}{L}$-co-monotone \cite{bauschke2020generalized} (also known as $-\frac{1}{L}$-cohypomonotone).

\textbf{Resolvent operator.}
The operator $J_Gx := \set{y \in \R^p : x \in y + Gy}$ is called the resolvent of $G$, often denoted by $J_Gx = (\Id + G)^{-1}x$, where $\Id$ is the identity mapping.
Clearly, evaluating $J_G$ requires solving a strongly monotone inclusion $0 \in y - x + Gy$ in $y$ for given $x$.
If $G$ is monotone, then $J_G$ is singled-valued, and if $G$ is maximally monotone, then $J_G$ is singled-valued and $\dom{J_G} = \R^p$.
If $G$ is monotone, then $J_G$ is firmly nonexpansive \cite[Proposition 23.10]{Bauschke2011}. 

\beforesec
\section{The Equivalence of Halpern's and Nesterov's Accelerated Schemes}\label{sec:main_part1}
\aftersec
Both Nesterov's accelerated and Halpern's fixed-point iteration schemes show significant improvement on convergence rates over classical methods for solving \eqref{eq:3o_MI}.
However, they are derived from different perspectives, and it is unclear if they are closely related to each other.
In this section, we show that these schemes are actually equivalent, though they may use different sets of parameters.

To present our analysis, we consider the following co-coercive equation: \hspace{-2ex}
\myeq{eq:MI}{
\text{Find}~y^{\star}\in\dom{G}~\text{such that:}~Gy^{\star} = 0,
}
where $G : \R^p\to\R^p$ is a single-valued and $\frac{1}{L}$-co-coercive operator.
We denote by $\zer{G} := G^{-1}(0) = \set{y^{\star} \in \dom{G} : Gy^{\star} = 0}$ the solution set of \eqref{eq:MI}, and assume that $\zer{G}$ is nonempty.

\beforesubsec
\subsection{The Halpern fixed-point scheme and its convergence}\label{subsec:HP_scheme00}
\aftersubsec
The Halpern fixed-point scheme \cite{diakonikolas2020halpern,halpern1967fixed,lieder2021convergence} for solving \eqref{eq:MI} is written as follows:
\myeq{eq:HP_scheme00}{
y_{k+1} := \beta_ky_0 + (1-\beta_k)y_k  - \eta_kGy_k,
}
where $\beta_k \in (0, 1)$ and $\eta_k > 0$ are appropriately chosen.

The convergence rate of \eqref{eq:HP_scheme00} has been established in \cite{diakonikolas2020halpern,lieder2021convergence} using different tools.
While \cite{lieder2021convergence} provides a direct proof and uses a performance estimation problem approach to establish convergence of \eqref{eq:HP_scheme00}, \cite{diakonikolas2020halpern} exploits a Lyapunov's technique to analyze its convergence rate.  
Let us summarize the result in \cite{diakonikolas2020halpern} in our context.

The standard Lyapunov function to study \eqref{eq:HP_scheme00} is 
\myeq{eq:HP_scheme00_Lyapunov_func}{
\Lc_k := \frac{p_k}{L}\norms{Gy_k}^2 + q_k\iprods{Gy_k, y_k - y_0},
}
where $p_k := q_0k(k+1)$ and $q_k := q_0(k+1)$ (for some $q_0 > 0$) are given parameters.
\revise{The following theorem is from \cite{diakonikolas2020halpern,lieder2021convergence} and states the convergence rate of \eqref{eq:HP_scheme00}.}

\begin{theorem}[\cite{diakonikolas2020halpern,lieder2021convergence}]\label{th:HP_scheme00_convergence}
Assume that $G$ in \eqref{eq:MI} is $\frac{1}{L}$-co-coercive with $L \in (0, +\infty)$, and $y^{\star} \in \zer{G}$.
Let $\sets{y_k}$ be generated by \eqref{eq:HP_scheme00} using $\beta_k := \frac{1}{k+2}$ and $\eta_k := \frac{2(1-\beta_k)}{L}$.
Then
\myeq{eq:HP_scheme00_convergence}{
\norms{Gy_k} \leq \frac{L\norms{y_0 - y^{\star}}}{k+1}.
}
\end{theorem}
\begin{remark}\label{re:th1_remark}
\revise{
If we choose $\beta_k := \frac{1}{k+2}$ and $\eta_k := \frac{1-\beta_k}{L}$, then using a similar proof as in Theorem~\ref{th:HP_scheme00_convergence} from \cite{diakonikolas2020halpern,lieder2021convergence}, we can show that
\myeqn{
\norms{Gy_k}^2 \leq \frac{4L^2\norms{y_0 - y^{\star}}^2}{(k+1)(k+3)} \ \ \text{and} \ \ \sum_{k=0}^{\infty}(k+1)(k+2)\norms{Gy_{k+1} - Gy_k}^2 \leq 2L^2\norms{y_0 - y^{\star}}^2.
}
However, if we choose $\eta_k := \frac{2(1-\beta_k)}{L}$ as in Theorem~\ref{th:HP_scheme00_convergence}, then we do not obtain the last summable inequality.
}
\end{remark}
Note that if $\beta_k := \frac{1}{k+2}$, then we can rewrite \eqref{eq:HP_scheme00} into the Halpern fixed-point iteration as in \cite{lieder2021convergence}:
\myeq{eq:HP_scheme00_b}{
y_{k+1} := \tfrac{1}{k+2}y_0 + \left(1 - \tfrac{1}{k+2}\right)Ty_k, \quad \text{where} \quad Ty_k := y_k - \tfrac{2}{L}Gy_k.
}
Since $G$ is $\frac{1}{L}$-co-coercive, $T = \Id - \frac{2}{L}G$ is nonexpansive, see \cite[Proposition 4.11]{Bauschke2011}.
Therefore, \eqref{eq:HP_scheme00} is equivalent to the scheme studied in \cite{lieder2021convergence}, and Theorem~\ref{th:HP_scheme00_convergence} can be obtained from the results in \cite{lieder2021convergence}.
The choice of $\beta_k$ and $\eta_k$ in Theorem~\ref{th:HP_scheme00_convergence} are tight and the bound \eqref{eq:HP_scheme00_convergence} is unimprovable since there exists an example that achieves this rate as the lower bound, see, e.g., \cite{lieder2021convergence} for such an example.

\beforesubsec
\subsection{The equivalence between Halpern's and Nesterov's accelerated schemes}\label{subsec:NAG2HP}
\aftersubsec
Our next step is to show that the Halpern fixed-point iteration \eqref{eq:HP_scheme00} can be transformed into a Nesterov's accelerated interpretation and vice versa.

\begin{theorem}\label{th:HP2Nes_scheme}
Let $\sets{x_k}$ and $\sets{y_k}$ be generated by the following scheme starting from $y_0 \in \R^p$ and $x_0 = x_{-1} = y_{-1} := y_0$ and $\beta_{-1} = \eta_{-1} = 0$:
\myeq{eq:NAG_scheme00_2corr}{
\left\{\begin{array}{lcl}
x_{k+1} &:= & y_k - \gamma_kGy_k, \vspace{1ex}\\
y_{k+1} &:= & x_{k+1} + \theta_k(x_{k+1} - x_k) + \nu_k(y_k - x_{k+1}) + \kappa_k(y_{k-1} - x_k),
\end{array}\right.
}
where $\theta_k := \frac{\beta_k(1-\beta_{k-1})}{\beta_{k-1}}$, $\nu_k := \frac{\beta_k}{\beta_{k-1}} + 1 - \beta_k - \frac{\eta_k}{\gamma_k}$, and $\kappa_k := \frac{\beta_k}{\beta_{k-1}}\left(\frac{\eta_{k-1}}{\gamma_{k-1}} -  1 + \beta_{k-1}\right)$.
Then, the sequence $\sets{y_k}$ is  identical  to the one generated by \eqref{eq:HP_scheme00} for solving \eqref{eq:MI}.

In particular, if  $\gamma_k := \frac{\eta_k}{1-\beta_k}$, then $\nu_k = \frac{\beta_k}{\beta_{k-1}}$, $\kappa_k = 0$, and \eqref{eq:NAG_scheme00_2corr} reduces to 
\myeq{eq:NAG_scheme00}{
\arraycolsep=0.3em
\left\{\begin{array}{lcl}
x_{k+1} &:= & y_k - \gamma_kGy_k,  \vspace{1ex}\\
y_{k+1} &:= & x_{k+1} + \theta_k(x_{k+1} - x_k) + \nu_k(y_k - x_{k+1}).
\end{array}\right.
}
\end{theorem}

\revise{Both \eqref{eq:NAG_scheme00_2corr} and \eqref{eq:NAG_scheme00} can be viewed as Nesterov's accelerated variants with correction terms.
Here, \eqref{eq:NAG_scheme00_2corr} has two correction terms $\nu_k(y_k - x_{k+1})$ and $\kappa_k(y_{k-1} - x_k)$, while \eqref{eq:NAG_scheme00} has only one term $\nu_k(y_k - x_{k+1})$.
In fact, \eqref{eq:NAG_scheme00_2corr} shows that  \eqref{eq:HP_scheme00} is equivalent to Nesterov's accelerated scheme with gradient correction in  \cite{shi2021understanding} as shown in Remark~\ref{re:HS2NAG_remark1} below.
Alternatively, compared to the ``optimized gradient method'' (OGM1) in \cite{kim2016optimized} for convex optimization, our coefficient $\nu_k$ in  \eqref{eq:NAG_scheme00} is positive in contrast to a negative value in OGM1.
Note that \eqref{eq:NAG_scheme00} covers the proximal-point scheme in \cite{mainge2021accelerated} as a special case.
As discussed in \cite{attouch2022ravine}, \eqref{eq:NAG_scheme00} can be viewed as a variant of Ravine's method if the convergence rate is given in $y_k$ instead of $x_k$.
}

\begin{proof}[\textbf{Proof of Theorem~\ref{th:HP2Nes_scheme}}]
\textbf{[$\eqref{eq:NAG_scheme00_2corr}\Rightarrow\eqref{eq:HP_scheme00}$]}
Substituting  $\theta_k$, $\nu_k$, and $\kappa_k$ into \eqref{eq:NAG_scheme00_2corr}, and simplifying the result, we get
\myeqn{ 
\arraycolsep=0.2em
\begin{array}{lcl}
y_{k+1} &= &\left(\frac{\beta_k}{\beta_{k-1}} - \beta_k + 1\right) x_{k+1} - \frac{\beta_k(1-\beta_{k-1})}{\beta_{k-1}}x_k +  \left( \frac{\beta_k}{\beta_{k-1}} + 1 - \beta_k - \frac{\eta_k}{\gamma_k}\right)(y_k - x_{k+1}) \vspace{1ex}\\
&& + {~} \frac{\beta_k}{\beta_{k-1}}\left(\frac{\eta_{k-1}}{\gamma_{k-1}} -  1 + \beta_{k-1}\right)(y_{k-1} - x_k).
\end{array}
}
Now, using the first line of \eqref{eq:NAG_scheme00_2corr} into this expression, we get
\myeqn{ 
\arraycolsep=0.2em
\begin{array}{lcl}
y_{k+1} &= & \left(\frac{\beta_k}{\beta_{k-1}} - \beta_k + 1\right)( y_k - \gamma_kGy_k) - \frac{\beta_k(1-\beta_{k-1})}{\beta_{k-1}}(y_{k-1} -  \gamma_{k-1}Gy_{k-1}) \vspace{1ex}\\
&& + {~} \left( \frac{\beta_k}{\beta_{k-1}} + 1 - \beta_k - \frac{\eta_k}{\gamma_k}\right)\gamma_kGy_k  + \frac{\beta_k}{\beta_{k-1}}\left(\frac{\eta_{k-1}}{\gamma_{k-1}} -  1 + \beta_{k-1}\right)\gamma_{k-1}Gy_{k-1}.
\end{array}
}
Rearranging this expression, we arrive at
\myeqn{
\arraycolsep=0.2em
\begin{array}{lcl}
\frac{1}{\beta_k}y_{k+1} - \left(\frac{1}{\beta_k} - 1\right)y_k + \frac{\eta_k}{\beta_k}Gy_k &= & \frac{1}{\beta_{k-1}}y_k - \left(\frac{1}{\beta_{k-1}} - 1\right)y_{k-1} +  \frac{\eta_{k-1}}{\beta_{k-1}}Gy_{k-1}.
\end{array}
}
By induction, and noticing that $y_{-1} = y_0$, and $\eta_{-1} = 0$, this expression leads to $\frac{1}{\beta_k}y_{k+1} - \big(\frac{1}{\beta_k} - 1\big)y_k + \frac{\eta_k}{\beta_k}Gy_k =  y_0$, which is indeed equivalent to \eqref{eq:HP_scheme00}.

\textbf{[\eqref{eq:HP_scheme00}$\Rightarrow$\eqref{eq:NAG_scheme00_2corr}]}
First, shifting the index from $k$ to $k-1$ in \eqref{eq:HP_scheme00}, we have $y_k = \beta_{k-1}y_0 + (1-\beta_{k-1})y_{k-1} - \eta_{k-1}Gy_{k-1}$.
Here, we assume that $y_{-1} = y_0$.
Multiplying this expression by $-\beta_k$ and \eqref{eq:HP_scheme00} by $\beta_{k-1}$ and adding the results, we obtain
\myeqn{
\beta_{k-1}y_{k+1} - \beta_ky_k = \beta_{k-1}(1-\beta_k)y_k - \beta_k(1-\beta_{k-1})y_{k-1} - \beta_{k-1}\eta_kGy_k + \beta_k\eta_{k-1}Gy_{k-1}.
}
This expression leads to
\myeq{eq:halpern_iter2}{
y_{k+1} = \left( \tfrac{\beta_k}{\beta_{k-1}} + 1 - \beta_k\right)y_k - \eta_kGy_k - \tfrac{\beta_k(1-\beta_{k-1})}{\beta_{k-1}}y_{k-1} + \tfrac{\beta_k\eta_{k-1}}{\beta_{k-1}}Gy_{k-1}.
}
Next, let us introduce $x_{k+1} := y_k - \gamma_kGy_k$ for some $\gamma_k > 0$.
Then, we have $Gy_k = \frac{1}{\gamma_k}(y_k - x_{k+1})$.
Substituting this relation into \eqref{eq:halpern_iter2}, we obtain 
\myeqn{
\arraycolsep=0.1em
\begin{array}{lcl}
y_{k+1} &= & \left( \frac{\beta_k}{\beta_{k-1}} + 1 - \beta_k\right)y_k - \frac{\eta_k}{\gamma_k}(y_k - x_{k+1}) - \frac{\beta_k(1-\beta_{k-1})}{\beta_{k-1}}y_{k-1} + \frac{\beta_k\eta_{k-1}}{\beta_{k-1}\gamma_{k-1}}(y_{k-1} - x_k) \vspace{1ex}\\
&= & x_{k+1} + \frac{\beta_k(1-\beta_{k-1})}{\beta_{k-1}}(x_{k+1} - x_k) +  \left( \frac{\beta_k}{\beta_{k-1}} + 1 - \beta_k - \frac{\eta_k}{\gamma_k}\right)(y_k - x_{k+1}) \vspace{1ex}\\
&& + {~} \frac{\beta_k}{\beta_{k-1}}\left(\frac{\eta_{k-1}}{\gamma_{k-1}} -  1 + \beta_{k-1}\right)(y_{k-1} - x_k).
\end{array}
}
If we let $\theta_k := \frac{\beta_k(1-\beta_{k-1})}{\beta_{k-1}}$, $\nu_k :=  \frac{\beta_k}{\beta_{k-1}} + 1 - \beta_k - \frac{\eta_k}{\gamma_k}$ and $\kappa_k := \frac{\beta_k}{\beta_{k-1}}\left(\frac{\eta_{k-1}}{\gamma_{k-1}} -  1 + \beta_{k-1}\right)$, then this expression can be rewritten as $y_{k+1} =  x_{k+1} + \theta_k(x_{k+1} - x_k) +  \nu_k(y_k - x_{k+1}) + \kappa_k(y_{k-1} - x_k)$.
Combining this line and $x_{k+1} = y_k - \gamma_kGy_k$, we get \eqref{eq:NAG_scheme00_2corr}.

Finally, if we choose $\gamma_k := \frac{\eta_k}{1 - \beta_k}$, then it is obvious that $\kappa_k = 0$, and $\nu_k = \frac{\beta_k}{\beta_{k-1}}$.
Hence, \eqref{eq:NAG_scheme00_2corr} reduces to \eqref{eq:NAG_scheme00}.
\Eproof
\end{proof}

\begin{remark}\label{re:HS2NAG_remark1}
\revise{Using \eqref{eq:halpern_iter2}, we can rewrite the Halpern-type scheme \eqref{eq:HP_scheme00} equivalently to
\myeqn{
y_{k+1} :=  y_k + \theta_k(y_k - y_{k-1}) - r_k Gy_k - s_k\left(Gy_k - Gy_{k-1}\right),
}
where $\theta_k :=  \frac{\beta_k(1-\beta_{k-1})}{\beta_{k-1}}$, $r_k := \eta_k - \tfrac{\beta_k\eta_{k-1}}{\beta_{k-1}}$, and $s_k := \tfrac{\beta_k\eta_{k-1}}{\beta_{k-1}} > 0$.
This expression shows that, under an appropriate choice of parameters, \eqref{eq:HP_scheme00} is equivalent to Nesterov's accelerated scheme with gradient correction studied in \cite{shi2021understanding}, when $Gy = \nabla{f}(y)$, the gradient of a convex function $f$.
In particular, if $\beta_k = \frac{1}{k+2}$ and $\eta_k := \frac{2(1-\beta_k)}{L}$ as in Theorem~\ref{th:HP_scheme00_convergence}, then $r_k$ reduces to $r_k = \frac{2}{(k+2)L} > 0$.
}
\end{remark}

\begin{remark}\label{re:HS2NAG_remark2}
If we choose $\gamma_k := \frac{1}{L}$, $\beta_k := \frac{1}{k+2}$, and $\eta_k := \frac{1 - \beta_k}{L} = \frac{k+1}{L(k+2)}$, then \eqref{eq:NAG_scheme00_2corr} reduces to the following one:
\myeq{eq:NAG_scheme00_a}{
\arraycolsep=0.3em
\left\{\begin{array}{lcl}
x_{k+1} &:= & y_k - \frac{1}{L}Gy_k, \vspace{1ex}\\
y_{k+1} &:= & x_{k+1} +  \frac{k}{k+2}(x_{k+1} - x_k) + \frac{k+1}{k+2}(y_k - x_{k+1}).  
\end{array}\right.
}
Alternatively, if $\gamma_k := \frac{1}{L}$, $\beta_k := \frac{1}{k+2}$, and $\eta_k := \frac{2(1 - \beta_k)}{L}$, then \eqref{eq:NAG_scheme00_2corr} becomes \hspace{-3ex}
\myeq{eq:NAG_scheme00_b}{
\arraycolsep=0.3em
\left\{\begin{array}{lcl}
x_{k+1} &:= & y_k - \frac{1}{L}Gy_k, \vspace{1ex}\\
y_{k+1} &:= & x_{k+1} +  \frac{k}{k+2}(x_{k+1} - x_k) + \frac{k}{k+2}(y_{k-1} - x_k).  
\end{array}\right.
}
The scheme \eqref{eq:NAG_scheme00_b} covers \cite{kim2021accelerated} as a special case when $Gy = J_{\lambda A}y$, the resolvent of a maximally monotone operator $\lambda A$, which is firmly nonexpansive.
After this work was published, we were recently informed that another related paper, \cite{contreras2023optimal} also points out the equivalence between the inertial scheme \eqref{eq:NAG_scheme00_b} and Halpern's iteration (\cite[Proposition 5]{contreras2023optimal}, but using a specific choice $\frac{1}{k+2}$ of the coefficient $\beta_k$.
\end{remark}

\revise{Utilizing Theorem~\ref{th:HP_scheme00_convergence} and Theorem~\ref{th:HP2Nes_scheme}, we obtain the following corollary showing the convergence rate of \eqref{eq:NAG_scheme00_2corr} under a particular choice of parameters. 
However, this choice of parameters is not necessarily unique as shown in Remark~\ref{re:th1_remark}.}


\begin{corollary}\label{co:HP_scheme00_convergence}
Assume that $G$ in \eqref{eq:MI} is $\frac{1}{L}$-co-coercive and $\zer{G}\neq\emptyset$.
Let $\set{(x_k, y_k)}$ be generated by \eqref{eq:NAG_scheme00_2corr} with $\gamma_k := \frac{1}{L}$, $\beta_k := \frac{1}{k+2}$, and $\eta_k := \frac{1-\beta_k}{L}$.
Then, we obtain $\theta_k := \frac{k}{k+2}$, $\nu_k := \frac{k+1}{k+2}$, and $\kappa_k := 0$.
Moreover, \eqref{eq:NAG_scheme00_2corr} reduces to \eqref{eq:NAG_scheme00} $($or equivalently \eqref{eq:NAG_scheme00_a}$)$, and the following guarantee holds:
\myeq{eq:HP_scheme00_convergence_b}{
\arraycolsep=0.3em
\begin{array}{rcl}
\norms{Gy_k}^2  & \leq & \frac{4L^2\norms{y_0 - y^{\star}}^2}{(k+1)(k+3)},  \vspace{1ex}\\
\sum_{k=0}^{\infty}(k+1)(k+2)\norms{Gy_{k+1} - Gy_k}^2 & \leq & 2L^2\norms{y_0 - y^{\star}}^2.
\end{array}
}
If we use $\gamma_k := \frac{1}{L}$, $\beta_k := \frac{1}{k+2}$, and $\eta_k := \frac{2(1-\beta_k)}{L}$, then we obtain $\theta_k := \frac{k}{k+2}$, $\nu_k := 0$, and $\kappa_k := \frac{k}{k+2}$, and  \eqref{eq:NAG_scheme00_2corr} reduces to  \eqref{eq:NAG_scheme00_b}.
Moreover, one has $\norms{Gy_k} \leq \frac{L\norms{y_0 - y^{\star}}}{k+1}$. 
\end{corollary}

The constant factor in the bound \eqref{eq:HP_scheme00_convergence_b} is slightly worse than the one in $\norms{Gy_k} \leq \frac{L\norms{y_0 - y^{\star}}}{k+1}$. 
In fact, the latter one is exactly optimal since there exists an example showing it matches the lower bound complexity, see, e.g. \cite{diakonikolas2020halpern,lieder2021convergence}.

\beforesubsec
\subsection{Convergence analysis of Nesterov's accelerated scheme \eqref{eq:NAG_scheme00}}\label{subsec:NAG0_analysis}
\aftersubsec
Now, we provide a direct analysis of \eqref{eq:NAG_scheme00} without using Theorem~\ref{th:HP_scheme00_convergence}.
For  simplicity, we will analyze the convergence of \eqref{eq:NAG_scheme00} with only one correction term.
However, our analysis can be easily extended to  \eqref{eq:NAG_scheme00_2corr} when $\kappa_k\neq 0$ with some simple modifications. 

Our analysis relies on the following Lyapunov function:
\myeq{eq:NAG_scheme00_Lyapunov_func1}{
\Vc_k := a_k\norms{Gy_{k-1}}^2 + b_k\iprods{Gy_{k-1}, x_k - y_k} + \norms{x_k  + t_k(y_k - x_k) - y^{\star} }^2 + \mu \norms{x_k - y^{\star}}^2,
}
where \revise{$a_k$, $b_k$, $t_k > 0$, and $\mu \geq 0$ are given parameters, which will be determined later.}
This Lyapunov is slightly different from $\Lc_k$ defined by \eqref{eq:HP_scheme00_Lyapunov_func}, but it is closely related to standard Nesterov's potential function (see, e.g., \cite{attouch2016rate,dAspremont2021acceleration,shi2021understanding}).
To see the connection between $\Vc_k$ and $\Lc_k$, we prove the following lemma.

\begin{lemma}\label{le:from_L_to_V}
Let $\Lc_k$ be defined by  \eqref{eq:HP_scheme00_Lyapunov_func} and $\Vc_k$ be defined by \eqref{eq:NAG_scheme00_Lyapunov_func1}.
Assume that $a_{k+1} :=  \frac{4p_k^2}{Lq_k^2} + \frac{4p_k\eta_k}{Lq_k(1-\beta_k)}$, $b_{k+1} := \frac{4p_k}{Lq_k\beta_k}$, $t_{k+1} := \frac{1}{\beta_k}$, and $\gamma_k := \frac{\eta_k}{1-\beta_k}$.
Then, we have
\myeq{eq:from_L_to_V}{
\Vc_{k+1} = \frac{4p_k}{Lq_k^2}\Lc_k  + \norms{y_0 - y^{\star}}^2 + \mu \norms{x_{k+1} - y^{\star}}^2.
}
\end{lemma}

Before proving this lemma, we make the following remarks.
\begin{itemize}
\item \revise{Note that the proof of Theorem~\ref{th:HP_scheme00_convergence} uses $\Lc_k$ defined by \eqref{eq:HP_scheme00_Lyapunov_func} with $p_k = q_0k(k+1)$ and $q_k = q_0(k+1)$.
In this case, we have $\frac{4p_k}{Lq_k^2} = \frac{4k}{Lq_0(k+1)} \approx \frac{4}{Lq_0}$.}
Hence, we can show that $\Vc_{k+1} = \frac{4k}{Lq_0(k+1)}\Lc_k +  \norms{y_0 - y^{\star}}^2 + \mu \norms{x_{k+1} - y^{\star}}^2$.
\item If we choose $p_k= c q_k^2$ for some $c > 0$, then $\Vc_{k+1} = \frac{4c}{L}\Lc_k +  \norms{y_0 - y^{\star}}^2 + \mu \norms{x_{k+1} - y^{\star}}^2$.
Clearly, if $\mu = 0$, then $\Vc_{k+1} = \frac{4c}{L}\Lc_k +  \norms{y_0 - y^{\star}}^2$, showing that $\Vc_k$ is equivalent to $\Lc_k$.

\item The term $\mu\norms{x_k - y^{\star}}^2$ allows us to get the tail $\norms{x_{k+1} - x_k}^2$ in \eqref{eq:NAG_scheme00_key_est2} of Lemma~\ref{le:NAG_scheme00_choice_of_params}, which is a key to prove convergence in Theorem~\ref{th:NAG_scheme00_convergence2}, especially the summable results and the $\SmallO{1/k^2}$-convergence rates.
It remains unclear to us how to prove such a convergence rate without using the term $\mu\norms{x_k - y^{\star}}^2$.
\end{itemize}

\begin{proof}[Proof of Lemma~\ref{le:from_L_to_V}]
First, from $\Lc_k = \frac{p_k}{L}\norms{Gy_k}^2 + q_k\iprods{Gy_k, y_k - y_0}$ in \eqref{eq:HP_scheme00_Lyapunov_func}, we can write it as 
\myeqn{
\arraycolsep=0.3em
\begin{array}{lcl}
\Lc_k &= & \frac{p_k}{L}\norms{Gy_k - \frac{Lq_k}{2p_k}(y_0 - y^{\star})}^2 + q_k\iprods{Gy_k, y_k - y^{\star}} - \frac{Lq_k^2}{4p_k}\norms{y_0 - y^{\star}}^2.
\end{array}
}
Second, from \eqref{eq:HP_scheme00}, we have $y_0 - y^{\star} = \frac{1}{\beta_k}(y_{k+1} + \eta_kGy_k - (1-\beta_k)y_k) - y^{\star} = y_k - y^{\star} + \frac{1}{\beta_k}(y_{k+1} - y_k) + \frac{\eta_k}{\beta_k}Gy_k$.
Using $x_{k+1} = y_k - \frac{\eta_k}{1-\beta_k} Gy_k$, we have $y_k = x_{k+1} + \frac{\eta_k}{1-\beta_k} Gy_k$.
Combining these lines, we get
\myeqn{ 
\arraycolsep=0.3em
\begin{array}{lcl}
Gy_k - \frac{Lq_k}{2p_k}(y_0 - y^{\star}) = Gy_k -  \frac{Lq_k}{2p_k} \big[ x_{k+1}  + \frac{1}{\beta_k}(y_{k+1} - x_{k+1}) - y^{\star} \big].
\end{array}
}
Now, let $z_{k+1} := x_{k+1} +  \frac{1}{\beta_k}(y_{k+1} - x_{k+1})$.
\revise{Then, we can rewrite $\Lc_k$ as 
\myeqn{
\arraycolsep=0.3em
\begin{array}{lcl}
\Lc_k &= & \frac{p_k}{L}\norms{Gy_k - \frac{Lq_k}{2p_k}(z_{k+1} - y^{\star})}^2 + q_k\iprods{Gy_k, y_k - y^{\star}} - \frac{Lq_k^2}{4p_k}\norms{y_0 - y^{\star}}^2.
\end{array}
}
Since $y_k = x_{k+1} + \frac{\eta_k}{1-\beta_k} Gy_k$, using the definition of $z_{k+1}$, we have $y_k - z_{k+1} = x_{k+1} - z_{k+1} +  \frac{\eta_k}{1-\beta_k} Gy_k = \frac{1}{\beta_k}(x_{k+1} - y_{k+1}) + \frac{\eta_k}{1-\beta_k} Gy_k$.
Further expanding the first term of $\Lc_k$ and using this relation, we can easily show that
\myeqn{
\arraycolsep=0.2em
\begin{array}{lcl}
\Lc_k &= & \revise{ \frac{p_k}{L}\norms{Gy_k}^2 + q_k\iprods{Gy_k, y_k - z_{k+1}} + \frac{Lq_k^2}{4p_k}\norms{z_{k+1} - y^{\star}}^2 - \frac{Lq_k^2}{4p_k}\norms{y_0 - y^{\star}}^2 } \vspace{1ex}\\
&= & \frac{Lq_k^2}{4p_k}\left[ \left( \frac{4p_k^2}{L^2q_k^2} + \frac{4p_k\eta_k}{Lq_k(1-\beta_k)}\right)\norms{Gy_k}^2 + \frac{4p_k}{Lq_k\beta_k}\iprods{Gy_k, x_{k+1} - y_{k+1}} + \norms{z_{k+1} - y^{\star}}^2 \right] \vspace{1ex}\\
&&  - {~} \frac{Lq_k^2}{4p_k}\norms{y_0 - y^{\star}}^2.
\end{array}
}
}
Finally, this expression together with \eqref{eq:NAG_scheme00_Lyapunov_func1} imply \eqref{eq:from_L_to_V}.
\Eproof
\end{proof}

Now, we prove the following key lemma for our convergence analysis.

\begin{lemma}\label{le:NAG_scheme00_key_est}
Let $\sets{(x_k, y_k)}$ be generated by \eqref{eq:NAG_scheme00} using $\gamma_k := \gamma > 0$, and $\Vc_k$ be defined by \eqref{eq:NAG_scheme00_Lyapunov_func1}.
Then, if $b_{k+1}\theta_k +  2\gamma t_k(t_k - 1) - 2\gamma \nu_k \theta_kt_{k+1}^2   \geq 0$, we have
\myeq{eq:NAG_scheme00_key_est}{
\arraycolsep=0.1em
\hspace{-4ex}
\begin{array}{lcl}
\Vc_k - \Vc_{k+1} &\geq & \left( \gamma b_{k+1}\nu_k  + \gamma^2 t_k^2 -  \gamma^2 t_{k+1}^2\nu_k^2 - a_{k+1} - \frac{\gamma^2b_k^2}{4a_k} \right) \norms{Gy_k}^2 \vspace{1ex}\\
&&+ {~}  \left[ b_{k+1}\theta_k +  2\gamma t_k(t_k - 1) - 2\gamma\nu_k\theta_kt_{k+1}^2 - b_k \right] \iprods{Gy_{k-1}, x_{k+1} - x_k}   \vspace{1ex}\\
&& + {~}  \left( t_k^2 - 2t_k + 1 + \mu - t_{k+1}^2\theta_k^2 \right) \norms{x_{k+1} - x_k}^2  \vspace{1ex}\\
&& + {~} \left( \frac{1}{L} - \gamma \right) \left[ b_{k+1}\theta_k +  2\gamma  t_k(t_k - 1) - 2\gamma\nu_k \theta_k t_{k+1}^2  \right] \norms{ Gy_k - Gy_{k-1} }^2 \vspace{1ex}\\
&& + {~} 2(t_k - t_{k+1}\theta_{k} - 1 - \mu \big)\iprods{x_{k+1} - x_k, x_{k+1} - y^{\star}} \vspace{1ex}\\
&& +  {~}  2\gamma\left(t_k - t_{k+1}\nu_k \right) \iprods{Gy_k, x_{k+1} - y^{\star}} \vspace{1ex}\\
&& + {~} a_k\big\Vert Gy_{k-1} - \frac{\gamma b_k}{2a_k}Gy_k\big\Vert^2. 
\end{array}
\hspace{-3ex}
}
\end{lemma}

\begin{proof}
Firstly, from \eqref{eq:NAG_scheme00_Lyapunov_func1}, we have
\myeq{eq:Vk_Vk1}{
\hspace{-2ex}
\arraycolsep=0.2em
\begin{array}{lcl}
\Vc_k - \Vc_{k+1} &= & a_k\norms{Gy_{k-1}}^2 - a_{k+1}\norms{Gy_k}^2 + b_k\iprods{Gy_{k-1}, x_k - y_k}  \vspace{1ex}\\
&& - {~} b_{k+1}\iprods{Gy_k, x_{k+1} - y_{k+1}} + \mu \norms{x_k - y^{\star}}^2 - \mu \norms{x_{k+1} - y^{\star}}^2  \vspace{1ex}\\
&&  + {~}  \norms{x_k - y^{\star} + t_k(y_k - x_k)}^2 - \norms{x_{k+1} - y^{\star} + t_{k+1}(y_{k+1} - x_{k+1})}^2.
\end{array}
\hspace{-3ex}
}
Next, since $y_k = x_{k+1} + \gamma Gy_k$ from \eqref{eq:NAG_scheme00}, it is easy to show that
\myeq{eq:Vk_Vk1_proof1}{
\arraycolsep=0.2em
\begin{array}{lcl}
\iprods{Gy_{k-1}, x_k - y_k} &= & -\iprods{Gy_{k-1}, x_{k+1} - x_k} - \gamma \iprods{Gy_{k-1}, Gy_k}.
\end{array}
}
Similarly, from \eqref{eq:NAG_scheme00} we have $x_{k+1} - y_{k+1} = -\theta_k(x_{k+1} - x_k) - \gamma\nu_kGy_k$, leading to{\!\!\!}
\myeq{eq:Vk_Vk1_proof2}{
\iprods{Gy_k, x_{k+1} - y_{k+1}}  =  - \theta_k \iprods{Gy_k, x_{k+1} - x_k} -  \gamma\nu_k\norms{Gy_k}^2.
}
Then, using again $y_k = x_{k+1} + \gamma Gy_k$ from \eqref{eq:NAG_scheme00}, we can derive 
\myeq{eq:Vk_Vk1_proof3}{
\hspace{-2ex}
\arraycolsep=0.2em
\begin{array}{lcl}
\norms{x_k - y^{\star} + t_k(y_k - x_k)}^2 & = &  \norms{x_{k+1} - y^{\star} + (t_k-1)(x_{k+1} - x_k)  + \gamma t_k Gy_k}^2 \vspace{1ex}\\
&= & \norms{x_{k+1} - y^{\star}}^2 + (t_k - 1)^2 \norms{x_{k+1} - x_k}^2   \vspace{1ex}\\
&& + {~}  \gamma^2t_k^2\norms{Gy_k}^2 +  2(t_k-1)\iprods{x_{k+1} - x_k, x_{k+1} - y^{\star}}  \vspace{1ex}\\
&& + {~} 2\gamma t_k(t_k-1)\iprods{Gy_k, x_{k+1} - x_k} + 2\gamma t_k \iprods{Gy_k, x_{k+1} - y^{\star}}.
\end{array}
\hspace{-7ex}
}
Similarly, using $y_{k+1} - x_{k+1} = \theta_k(x_{k+1} - x_k) + \gamma\nu_kGy_k$, we can show that
\myeq{eq:Vk_Vk1_proof4}{
\hspace{-3ex}
\arraycolsep=0.2em
\begin{array}{lcl}
\norms{x_{k+1} - y^{\star} + t_{k+1}(y_{k+1} - x_{k+1})}^2 &= & \norms{x_{k+1} - y^{\star} + t_{k+1}\theta_k(x_{k+1} - x_k) + \gamma t_{k+1}\nu_kGy_k}^2 \vspace{1ex}\\
&= & \norms{x_{k+1} - y^{\star}}^2  + t_{k+1}^2\theta_k^2\norms{x_{k+1} - x_k}^2   \vspace{1ex}\\
&& + {~}  2\gamma\nu_k\theta_kt_{k+1}^2\iprods{Gy_k, x_{k+1} - x_k} \vspace{1ex}\\
&& + {~} 2t_{k+1}\theta_k \iprods{x_{k+1} - x_k, x_{k+1} - y^{\star}} \vspace{1ex}\\
&& + {~} 2\gamma t_{k+1}\nu_k \iprods{Gy_k, x_{k+1} - y^{\star}} + \gamma^2t_{k+1}^2\nu_k^2\norms{Gy_k}^2.
\end{array}
\hspace{-10ex}
}
Substituting \eqref{eq:Vk_Vk1_proof1},  \eqref{eq:Vk_Vk1_proof2},  \eqref{eq:Vk_Vk1_proof3}, and  \eqref{eq:Vk_Vk1_proof4} into \eqref{eq:Vk_Vk1}, and using $\mu \norms{x_k - y^{\star}}^2 -  \mu \norms{x_{k+1} - y^{\star}}^2 = \mu  \norms{x_{k+1} - x_k}^2 - 2\mu \iprods{x_{k+1} - x_k, x_{k+1} - y^{\star}}$, we can show that
\myeq{eq:Vk_Vk1_proof5}{
\arraycolsep=0.2em
\begin{array}{lcl}
\Vc_k - \Vc_{k+1} &= & a_k\norms{Gy_{k-1}}^2 + \left[\gamma b_{k+1}\nu_k  + \gamma^2t_k^2 -  \gamma^2 t_{k+1}^2\nu_k^2 - a_{k+1} \right] \norms{Gy_k}^2  \vspace{1ex}\\
&& - {~} \gamma b_k\iprods{Gy_{k-1}, Gy_k} - b_k\iprods{Gy_{k-1}, x_{k+1} - x_k} \vspace{1ex}\\
&& + {~}  \left[ b_{k+1}\theta_k +  2\gamma t_k(t_k - 1) - 2\gamma\nu_k\theta_kt_{k+1}^2 \right] \iprods{Gy_k, x_{k+1} - x_k}   \vspace{1ex}\\
&& + {~}  \left[ \mu +  (t_k - 1)^2 - t_{k+1}^2\theta_k^2 \right] \norms{x_{k+1} - x_k}^2 \vspace{1ex}\\
&& + {~} 2\left(t_k - t_{k+1}\theta_k - 1 - \mu \right)\iprods{x_{k+1} - x_k, x_{k+1} - y^{\star}} \vspace{1ex}\\
&& + {~} 2\gamma\big(t_k - t_{k+1}\nu_k \big)\iprods{Gy_k, x_{k+1} - y^{\star}}.
\end{array}
}
By the $\frac{1}{L}$-co-coerciveness of $G$ and $x_{k+1} = y_k - \gamma Gy_k$ from \eqref{eq:NAG_scheme00}, we can derive 
\myeqn{ 
\arraycolsep=0.2em
\begin{array}{lcl}
\iprods{Gy_k - Gy_{k-1}, x_{k+1} - x_k} & \geq & \left( \frac{1}{L} - \gamma \right) \norms{ Gy_k - Gy_{k-1} }^2.
\end{array}
}
This inequality implies that $\iprods{Gy_k, x_{k+1} - x_k} \geq \iprods{Gy_{k-1}, x_{k+1} - x_k} + \left( \frac{1}{L} - \gamma \right) \norms{ Gy_k - Gy_{k-1} }^2$.
Finally, if we assume that $b_{k+1}\theta_k +  2\gamma t_k(t_k - 1) - 2\gamma \nu_k\theta_kt_{k+1}^2  \geq 0$, then by substituting the last inequality into \eqref{eq:Vk_Vk1_proof5} and using $a_k\norms{Gy_{k-1}}^2 - \gamma b_k\iprods{Gy_k, Gy_{k-1}} = a_k\norms{Gy_{k-1} - \frac{\gamma b_k}{2a_k}Gy_k}^2 - \frac{\gamma^2b_k^2}{4a_k}\norms{Gy_k}^2$, we obtain \eqref{eq:NAG_scheme00_key_est}.
\Eproof
\end{proof}

Our next lemma is to provide a choice of parameters such that $\Vc_k - \Vc_{k+1} \geq 0$.
\begin{lemma}\label{le:NAG_scheme00_choice_of_params}
Let $0 < \gamma \leq \frac{2}{L}$, $\mu \geq 0$, and $\omega \geq 1$ be given.
Let $\sets{(x_k, y_k)}$ be generated by \eqref{eq:NAG_scheme00} and $\Vc_k$ be defined by \eqref{eq:NAG_scheme00_Lyapunov_func1}.
Assume that $t_k$, $\theta_k$, $\nu_k$, $a_k$, and $b_k$ in \eqref{eq:NAG_scheme00} and \eqref{eq:NAG_scheme00_Lyapunov_func1} are chosen as follows:
\myeq{eq:NAG_scheme00_choice_of_params}{
\arraycolsep=0.2em
\begin{array}{ll}
& t_k := \tfrac{k + 2\omega + 1}{\omega}, \quad  \theta_k := \tfrac{t_k  - 1 - \mu}{t_{k + 1}}, \quad \nu_k := 1 - \tfrac{1}{t_{k+1}}, \vspace{1ex}\\
& b_k := 2\gamma t_k(t_k - 1), \quad\text{and}\quad a_k := \gamma^2t_k(t_k-1).
\end{array}
}
Then, it holds that
\myeq{eq:NAG_scheme00_key_est2}{
\hspace{-1ex}
\arraycolsep=0.2em
\begin{array}{lcl}
\Vc_k - \Vc_{k+1} &\geq & \mu(2t_k - 1 - \mu) \norms{x_{k+1} - x_k}^2  + \revise{\frac{\gamma(\omega - 1)}{\omega}\left(\frac{2}{L} - \gamma \right)} \norms{Gy_k}^2  \vspace{1ex}\\
&& + {~} \gamma\big(\frac{2}{L} - \gamma\big)t_k(t_k-1)\norms{Gy_k -  Gy_{k-1}}^2 \geq 0.
\end{array}
\hspace{-2ex}
}
Moreover, we have \revise{$\Vc_k \geq  \frac{\gamma(2-\gamma L)(t_k-1)}{L}\norms{Gy_{k-1}}^2 + \mu\norms{x_k - y^{\star}}^2  \geq 0$} and
\myeq{eq:NAG_scheme00_key_est3}{
\arraycolsep=0.2em
\left\{\begin{array}{lcl}
\sum_{k=0}^{\infty} \mu(2t_k - 1 - \mu)\norms{x_{k+1} - x_k}^2 & \leq & \Vc_0, \vspace{1ex}\\
\revise{\frac{\gamma(2-L\gamma)(\omega - 1)}{L\omega} } \sum_{k=0}^{\infty}\norms{Gy_k}^2 &\leq & \Vc_0, \vspace{1ex}\\
\revise{\frac{\gamma(2 - L\gamma)}{L}} \sum_{k=0}^{\infty} t_k(t_k-1) \norms{Gy_k - Gy_{k-1}}^2 & \leq & \Vc_0, \vspace{1ex}\\
\revise{\frac{(\omega+1)(2-L\gamma)}{L\gamma(2\omega-1)}\sum_{k=0}^{\infty} t_k^2 \norms{x_{k+1} - x_k - \theta_{k-1}(x_k - x_{k-1})}^2} & \leq & \Vc_0.
\end{array}\right.
}
\end{lemma}

\begin{proof}
Let us show that $t_k$, $\theta_k$, $\nu_k$, and $b_k$ chosen by \eqref{eq:NAG_scheme00_choice_of_params} satisfy 
\myeq{eq:NAG_scheme00_params_cond}{
t_k - t_{k+1}\theta_k - 1 - \mu =  0 \quad \text{and}\quad b_{k+1}\theta_k +  2\gamma t_k(t_k - 1) - 2\gamma \nu_k \theta_k t_{k+1}^2  - b_k =  0.
}
First, since $\theta_k = \frac{t_k - 1 - \mu}{t_{k+1}}$,  the first condition of \eqref{eq:NAG_scheme00_params_cond} holds.
Next, using $b_k = 2\gamma t_k(t_k - 1)$ and $\nu_k = 1 - \frac{1}{t_{k+1}}$, we can easily verify the second condition of \eqref{eq:NAG_scheme00_params_cond}.

Now, using \eqref{eq:NAG_scheme00_choice_of_params}, we can directly compute the following coefficients of \eqref{eq:NAG_scheme00_key_est}:
\myeqn{
\arraycolsep=0.2em
\left\{\begin{array}{lcl}
t_k^2 - 2t_k + 1 + \mu - t_{k+1}^2\theta_k^2 & = &  \mu(2t_k - 1 - \mu), \vspace{1ex}\\
t_k - t_{k+1}\nu_k & = &  \frac{\omega - 1}{\omega}, \vspace{1ex}\\
 b_{k+1}\theta_k +  2\gamma  t_k(t_k - 1) - 2\gamma\nu_k \theta_k t_{k+1}^2 &= & 2\gamma t_k(t_k-1),\vspace{1ex}\\
\gamma b_{k+1} \nu_k  + \gamma^2 t_k^2 -  \gamma^2 t_{k+1}^2\nu_k^2  - a_{k+1} - \frac{\gamma^2b_k^2}{4a_k} & = &  \frac{\gamma^2(\omega-1)}{\omega}.
\end{array}\right.
}
\revise{Substituting \eqref{eq:NAG_scheme00_params_cond} and these expressions into \eqref{eq:NAG_scheme00_key_est}, and then using $a_k = \gamma^2t_k(t_k-1)$ and $2a_k = \gamma b_k$ from \eqref{eq:NAG_scheme00_choice_of_params} into the resulting inequality, we can simplify it as
\myeqn{
\arraycolsep=0.1em
\hspace{-4ex}
\begin{array}{lcl}
\Vc_k - \Vc_{k+1} &\geq &  \frac{\gamma^2(\omega-1)}{\omega} \norms{Gy_k}^2 +   \mu(2t_k - 1 - \mu) \norms{x_{k+1} - x_k}^2  \vspace{1ex}\\
&& + {~}  \gamma\left( \frac{2}{L} - \gamma \right)t_k(t_k-1)  \norms{ Gy_k - Gy_{k-1} }^2 + \frac{2\gamma(\omega - 1)}{\omega}  \iprods{Gy_k, x_{k+1} - y^{\star}}.
\end{array}
\hspace{-3ex}
}}
Moreover, since $x_{k+1} = y_k - \gamma Gy_k$, using $Gy^{\star} = 0$ and the $\frac{1}{L}$-co-coerciveness of $G$, we have $\iprods{Gy_k, x_{k+1} - y^{\star}} = \iprods{Gy_k, y_k - y^{\star}} - \gamma\norms{Gy_k}^2  \geq \left(\frac{1}{L} - \gamma\right)\norms{Gy_k}^2$.
\revise{Substituting this inequality and \eqref{eq:NAG_scheme00_params_cond} into the last expression, we get \eqref{eq:NAG_scheme00_key_est2}.}

\revise{Next, since $Gy^{\star} = 0$, using the $\frac{1}{L}$-co-coerciveness of $G$, we get 
\myeqn{
\arraycolsep=0.1em
\begin{array}{lcl}
\iprods{Gy_{k-1}, x_k - y^{\star}} & = & \iprods{Gy_{k-1} - Gy^{\star}, y_{k-1} - y^{\star} - \gamma (Gy_{k-1} - Gy^{\star})} \vspace{1ex}\\
& \geq & \left(\frac{1}{L} - \gamma\right)\norms{Gy_{k-1}}^2.
\end{array}
}
Utilizing this bound and the Cauchy-Schwarz inequality, we can show that
\myeqn{
\hspace{-2ex}
\arraycolsep=0.2em
\begin{array}{lcl}
b_k\iprods{Gy_{k-1}, x_k - y_k}  &= & \frac{b_k}{t_k}\iprods{Gy_{k-1}, x_k - y^{\star}}  - \frac{b_k}{t_k}\iprods{Gy_{k-1}, x_k - y^{\star} + t_k(y_k - x_k)} \vspace{1ex}\\
& \geq & -\frac{b_k^2}{4t_k^2}\norms{Gy_{k-1}}^2 - \norms{x_k - y^{\star} + t_k(y_k - x_k)}^2 \vspace{1ex}\\
&& + {~} \frac{b_k}{t_k}\left(\frac{1}{L}  -  \gamma\right)\norms{Gy_{k-1}}^2.
\end{array}
\hspace{-2ex}
}}
Substituting this bound into the definition \eqref{eq:NAG_scheme00_Lyapunov_func1} of $\Vc_k$ and noticing that $a_k - \frac{b_k^2}{4t_k^2} = \gamma^2(t_k - 1)$, we get $\Vc_k \geq \left(a_k - \frac{b_k^2}{4t_k^2}\right)\norms{Gy_{k-1}}^2 +  \frac{b_k}{t_k}\left(\frac{1}{L} - \gamma\right)\norms{Gy_{k-1}}^2 + \mu\norms{x_k - y^{\star}}^2  =   \frac{\gamma(2-\gamma L)(t_k-1)}{L}\norms{Gy_{k-1}}^2 + \mu\norms{x_k - y^{\star}}^2$, which proves that $\Vc_k \geq 0$.

\revise{Summing up \eqref{eq:NAG_scheme00_key_est2} from $k := 0$ to $k := K$ and using $\Vc_{K+1} \geq 0$, we get
\myeq{eq:NAG_scheme00_proof5}{
\hspace{-1ex}
\arraycolsep=0.2em
\begin{array}{ll}
\sum_{k=0}^K\big[\!\! & \mu(2t_k - 1 - \mu) \norms{x_{k+1} - x_k}^2  +  \frac{\gamma(\omega - 1)}{\omega}\left(\frac{2}{L} - \gamma \right) \norms{Gy_k}^2  \vspace{1ex}\\
& + {~} \gamma\big(\frac{2}{L} - \gamma\big)t_k(t_k-1)\norms{Gy_k -  Gy_{k-1}}^2 \big] \leq \Vc_0 - \Vc_{K+1} \leq \Vc_0.
\end{array}
\hspace{-2ex}
}
Letting $K \to \infty$ in this inequality, we obtain the first three expressions  of \eqref{eq:NAG_scheme00_key_est3}.
}

\revise{Finally, since $x_{k+1} - x_k - \theta_{k-1}(x_k - x_{k-1}) = \gamma(Gy_k - \nu_{k-1}Gy_{k-1})$, using Young's inequality and then $c_k := \frac{\omega t_k}{(\omega-1)(t_k-1)}$, we can derive that
\myeqn{
\arraycolsep=0.2em
\begin{array}{ll}
& \norms{x_{k+1} - x_k - \theta_{k-1}(x_k - x_{k-1})}^2   =   \gamma^2\norms{Gy_k - \nu_{k-1}Gy_{k-1}}^2 \vspace{1ex}\\
& \quad\quad \leq   \gamma^2(1+c_k)\big[ \nu_{k-1}^2 \norms{Gy_k - Gy_{k-1}}^2  +  \frac{(\nu_{k-1}-1)^2 }{c_k}\norms{Gy_k}^2 \big] \vspace{1ex}\\
& \quad\quad \leq   \gamma^2(1+c_k)\big[ \frac{(t_k-1)^2}{t_k^2} \norms{Gy_k - Gy_{k-1}}^2  +  \frac{1}{c_kt_k^2}\norms{Gy_k}^2 \big] \vspace{1ex}\\
& \quad\quad =   \frac{\gamma^2(1+c_k)(t_k-1)}{t_k^3}\big[ (t_k-1)t_k  \norms{Gy_k - Gy_{k-1}}^2  +  \frac{t_k}{c_k(t_k-1)}\norms{Gy_k}^2 \big] \vspace{1ex}\\
& \quad\quad \leq  \frac{\gamma^2(2\omega-1)}{(\omega-1)t_k^2}\big[ t_k(t_k-1) \norms{Gy_k - Gy_{k-1}}^2  + \frac{\omega-1}{\omega}\norms{Gy_k}^2 \big].
\end{array}
}
Combining this inequality and \eqref{eq:NAG_scheme00_proof5} (after dropping its first term), and then letting $K\to\infty$, we obtain the last line of \eqref{eq:NAG_scheme00_key_est3}.}
\Eproof
\end{proof}

The following theorem  proves  the convergence of Nesterov's accelerated scheme \eqref{eq:NAG_scheme00}, but using a different set of parameters compared to Theorem~\ref{th:HP_scheme00_convergence}.

\begin{theorem}\label{th:NAG_scheme00_convergence2}
Assume that $G$ in \eqref{eq:MI} is $\frac{1}{L}$-co-coercive and $\zer{G}\neq\emptyset$.
Let $\set{(x_k, y_k)}$ be generated by \eqref{eq:NAG_scheme00} to solve \eqref{eq:MI} using $\gamma_k := \gamma \in \left(0, \frac{1}{L}\right)$, $\theta_k := \frac{k+1}{k+ 2\omega + 2}$, and $\nu_k :=  \frac{k + \omega + 2}{k+ 2\omega + 2} \in (0, 1)$ for a given constant $\omega > 2$.
Then, we have
\myeq{eq:NAG_scheme00_convergence2}{
\arraycolsep=0.2em
\left\{\begin{array}{lcl}
\sum_{k=0}^{\infty}(k + \omega + 1)\norms{x_{k+1} - x_k}^2 < +\infty \quad \text{and} \quad \norms{x_{k+1} - x_k}^2 = \SmallO{\frac{1}{k^2}}, \vspace{1ex}\\
\sum_{k=0}^{\infty}(k + 2\omega + 1)\norms{y_k - x_k}^2 < +\infty \quad \text{and} \quad \norms{y_k - x_k}^2 = \SmallO{\frac{1}{k^2}}, \vspace{1ex}\\
\sum_{k=0}^{\infty}(k + \omega + 1)\norms{Gy_{k}}^2 < +\infty \quad \text{and} \quad \norms{Gy_k}^2 = \SmallO{\frac{1}{k^2}}, \vspace{1ex}\\
\sum_{k=0}^{\infty}(k + \omega + 1)\norms{Gx_k}^2 < +\infty \quad \text{and} \quad \norms{Gx_k}^2 = \SmallO{\frac{1}{k^2}}, \vspace{1ex}\\
\sum_{k=0}^{\infty}(k+\omega)\norms{y_{k+1} - y_k}^2 < +\infty \quad \text{and} \quad \norms{y_{k+1}  - y_k}^2 = \SmallO{\frac{1}{k^2}}.
\end{array}\right.
}
Consequently, both $\set{x_k}$ and $\set{y_k}$ converge to $y^{\star} \in \zer{G}$.
\end{theorem}

Before proving Theorem~\ref{th:NAG_scheme00_convergence2}, we make the following remarks.
\begin{remark}\label{re:th3_remark1}
\revise{First, if we choose $\gamma = \frac{2}{L}$, then we only obtain the first result of \eqref{eq:NAG_scheme00_convergence2} and $\norms{x_{k+1} - x_k}^2 = \SmallO{\frac{1}{k^2}}$.
This rate is theoretically better than the $\BigO{1/k^2}$ rate in \cite{kim2021accelerated} when $k$ is sufficiently large.
Second, although we only state the $\SmallO{\cdot}$ rates of the four different quantities in Theorem~\ref{th:NAG_scheme00_convergence2}, the corresponding $\BigO{\cdot}$ rates of these quantities can also be achieved through our proof below.
Moreover, the upper bound of these rates can be expressed explicitly.
For instance, by the first line of \eqref{eq:NAG_scheme00_key_est3} and \eqref{eq:NAG_scheme00_th1_proof1} below, we can easily prove a $\BigO{1/k^2}$-rate of $\norms{y_k-x_k}^2$ as
\begin{equation*}
\norms{y_k - x_k}^2 \leq \frac{(\omega+1)^2\Vc_0}{(k+2\omega+1)^2}, \quad\text{where $\Vc_0$ is given by \eqref{eq:NAG_scheme00_Lyapunov_func1}}.
\end{equation*}
Third, if $\gamma \in \left(0, \frac{1}{L}\right)$, then we can prove $\SmallO{\frac{1}{k^2}}$ convergence  rates of $\norms{Gy_k}^2$, $\norms{Gx_k}^2$, $\norms{y_k - x_k}^2$, and $\norms{y_{k+1} - y_k}^2$.
Finally, note that we can simply choose $\omega = 3$ to further simplify the results.
In this case, we obtain $\theta_k = \frac{k+1}{k+8}$, which is different from $\theta_k = \frac{k}{k+2}$ in  \eqref{eq:NAG_scheme00_a} obtained by Theorem~\ref{th:HP_scheme00_convergence}.
}
\end{remark}

\revise{The key step to prove Theorem~\ref{th:NAG_scheme00_convergence2} is Lemma~\ref{le:NAG_scheme00_key_est}.
We believe that this lemma is new and its proof is relatively elementary.
This proof technique can be further extended to study other methods.
For instance, it has been recently exploited to study accelerated randomized coordinate methods for solving \eqref{eq:MI} in \cite{tran2022accelerated}.
Moreover, it is worth to emphasize that our results in this paper (see, e.g., Corollary \ref{co:HP_convergence2} below) show that  Halpern's fixed-point methods for solving \eqref{eq:MI} can achieve both $\BigO{\cdot}$ and $\SmallO{\cdot}$ rates by choosing different parameters.
%
Note that $\SmallO{\cdot}$ convergence rates for Nesterov's accelerated methods have been recently studied in a number of works such as \cite{attouch2020convergence,attouch2016rate,mainge2021accelerated,mainge2021fast,labarre2022first}, but we are not aware of any $\SmallO{1/k^2}$ rate for Halpern's fixed-point methods a prior to this work.
After the first version of this paper is completed and available online, we find a recent work \cite{bot2022fast} that also studies $o(\cdot)$ convergence rates of different Nesterov's accelerated schemes for \eqref{eq:MI}.
However,  \cite{bot2022fast} relies on the discretizations of dynamical systems as in \cite{attouch2020convergence}, and thus is different from our approach.
Another related work is \cite{yoon2022accelerated}, which shows that the iterate sequences generated by several Halpern-type schemes are actually close to each other and eventually converge to a solution of \eqref{eq:MI}.
Though this work did not specifically study $o(\cdot)$, but several summable bounds were also obtained.
}

\begin{proof}[The proof of Theorem~\ref{th:NAG_scheme00_convergence2}]
\revise{For simplicity of our analysis, we fix $\mu := 1$ in $\Vc_k$ defined by \eqref{eq:NAG_scheme00_Lyapunov_func1}.}
The first claim in the first line of \eqref{eq:NAG_scheme00_convergence2} comes directly from \eqref{eq:NAG_scheme00_key_est3} by noticing that $t_k - 1 = \frac{k+\omega+1}{\omega}$.
Now, we prove the second line of \eqref{eq:NAG_scheme00_convergence2}.
Indeed, from \eqref{eq:NAG_scheme00}, we have
\myeqn{
\arraycolsep=0.2em
\begin{array}{lcl}
y_{k+1} - x_{k+1} & = &  \nu_k(x_{k+1} + \gamma Gy_k - x_k) + (\theta_k - \nu_k)(x_{k+1} - x_k) \vspace{1ex}\\
& = & \nu_k(y_k - x_k) + (\theta_k - \nu_k)(x_{k+1} - x_k).
\end{array}
}
Hence, by Young's inequality, $\nu_k \in (0, 1)$, and this expression, we can show that
\myeqn{
\arraycolsep=0.2em
\begin{array}{lcl}
t_{k+1}^2\norms{y_{k+1} - x_{k+1}}^2 & = &  t_{k+1}^2\norms{\nu_k(y_k - x_k) + (\theta_k - \nu_k)(x_{k+1} - x_k)}^2 \vspace{1ex}\\
&\leq & t_{k+1}^2\nu_k\norms{y_k - x_k}^2 + \frac{t_{k+1}^2(\theta_k - \nu_k)^2}{1-\nu_k}\norms{x_{k+1} - x_k}^2.
\end{array}
}
Notice from \eqref{eq:NAG_scheme00_choice_of_params} that $t_{k+1}^2\nu_k = t_k^2 - \frac{ \omega(\omega-2)t_k + \omega - 1}{\omega^2}$ 
and $\frac{t_{k+1}^2(\theta_k - \nu_k)^2}{1-\nu_k} = \frac{(\omega + 1)^2(k+2\omega + 2)}{\omega^3}$.
Utilizing these expressions into the last inequality, we obtain
\myeq{eq:NAG_scheme00_th1_proof1}{
\hspace{-1ex}
\arraycolsep=0.2em
\begin{array}{lcl}
\frac{\omega(\omega-2)t_k + \omega - 1}{\omega^2} \norms{y_k - x_k}^2 & \leq & t_k^2\norms{y_k - x_k}^2 - t_{k+1}^2\norms{y_{k+1} - x_{k+1}}^2 \vspace{1ex}\\
&& + {~} \frac{(\omega + 1)^2t_{k+1}}{\omega^2}\norms{x_{k+1} - x_k}^2.
\end{array}
\hspace{-3ex}
}
Summing up this estimate from $k := 0$ to $k := K$, we get
\myeqn{
\hspace{-2ex}
\arraycolsep=0.2em
\begin{array}{lcl}
\sum_{k=0}^{K}\frac{\omega(\omega-2)t_k + \omega - 1}{\omega^2}\norms{y_k - x_k}^2 & \leq &   \frac{(\omega + 1)^2}{\omega^3}\sum_{k=0}^K (k+2\omega+2)\norms{x_{k+1} - x_k}^2 \vspace{1ex}\\
&& + {~} t_0^2\norms{y_0 - x_0}^2.
\end{array}
\hspace{-3ex}
}
Using the first line of \eqref{eq:NAG_scheme00_convergence2} into this inequality and $\omega > 2$, we obtain $\sum_{k=0}^{\infty}\big[(\omega-2)(k+2\omega+1) + \omega - 1\big] \norms{y_k - x_k}^2 < +\infty$, which implies the first claim in the second line of \eqref{eq:NAG_scheme00_convergence2}.
Moreover, \eqref{eq:NAG_scheme00_th1_proof1} also shows that $\lim_{k\to\infty}t_k^2\norms{x_k - y_k}^2$ exists (Here, we use \cite[Lemma 2.5]{xu2002iterative}).
Combining this fact and $\sum_{k=0}^{\infty}(k + 2\omega + 1)\norms{y_k - x_k}^2 < +\infty$, we obtain $\lim_{k\to\infty}t_k^2\norms{x_k - y_k}^2 = 0$, which shows that $\norms{x_k - y_k}^2 = \SmallO{1/k^2}$.

To prove the third line of  \eqref{eq:NAG_scheme00_convergence2}, we note that $\gamma\nu_kGy_k = (y_{k+1} - x_{k+1}) - \theta_k(x_{k+1} - x_{k})$.
Hence, we have  $\gamma^2\nu_k^2(t_k - 1)\norms{Gy_k}^2 \leq 2(t_k-1)\norms{y_{k+1} - x_{k+1}}^2 + 2\theta_k^2(t_k-1)\norms{x_{k+1} - x_{k}}^2$.
Exploiting the last two terms from \eqref{eq:NAG_scheme00_convergence2}, we obtain $\sum_{k=0}^{\infty}(k + \omega + 1)\norms{Gy_k}^2 < + \infty$.

\revise{To prove the second part in the first line of  \eqref{eq:NAG_scheme00_convergence2},  utilizing both lines of \eqref{eq:NAG_scheme00}, we have $\theta_{k-1}^2\norms{x_k - x_{k-1}}^2 = \norms{y_k - x_k - \nu_{k-1}(y_{k-1} - x_k)}^2 =  \norms{x_{k+1} - x_k + \gamma(Gy_k - \nu_{k-1}Gy_{k-1})}^2 $.
Therefore, expanding this expression and using  $\iprods{Gy_k - Gy_{k-1}, x_{k+1} - x_k}  \geq  \big( \frac{1}{L} - \gamma \big) \norms{ Gy_k - Gy_{k-1} }^2$ from the proof of Lemma~\ref{le:NAG_scheme00_key_est}, we can derive that
\myeqn{
\hspace{-2ex}
\arraycolsep=0.1em
\begin{array}{lcl}
\Tc_{[1]} &: = & \theta_{k-1}^2t_k^2\norms{x_k - x_{k-1}}^2 - \theta_k^2t_{k+1}^2\norms{x_{k+1} - x_k}^2 \vspace{1ex}\\
&= &   \gamma^2 t_k^2\norms{Gy_k - \nu_{k-1}Gy_{k-1}}^2 + 2\gamma t_k^2\iprods{Gy_k - Gy_{k-1}, x_{k+1} - x_k} \vspace{1ex}\\
&& + {~} 2\gamma t_k^2(1-\nu_k)\iprods{Gy_{k-1}, x_{k+1} - x_k} + (t_k^2 - \theta_k^2t_{k+1}^2)\norms{x_{k+1} - x_k}^2 \vspace{1ex}\\
&\geq & \gamma^2 t_k^2\norms{Gy_k - \nu_{k-1}Gy_{k-1}}^2 + 2\gamma t_k^2\left(\frac{1}{L} - \gamma\right)\norms{Gy_k - Gy_{k-1}}^2 \vspace{1ex}\\
&& + {~} 2\gamma t_k^2(1-\nu_k)\iprods{Gy_{k-1}, x_{k+1} - x_k} + (t_k^2 - \theta_k^2t_{k+1}^2)\norms{x_{k+1} - x_k}^2.
\end{array}
\hspace{-3ex}
}
Note that by the choice of $t_k$, $\theta_k$, and $\nu_k$ as in Theorem~\ref{th:NAG_scheme00_convergence2}, we have $2\gamma t_k^2(1-\nu_k) \geq 0$ and $t_k^2 - \theta_k^2t_{k+1}^2 \geq 0$.
}
Employing the update rule \eqref{eq:NAG_scheme00_choice_of_params} and Young's inequality, the last inequality leads to 
\myeqn{
\arraycolsep=0.2em
\begin{array}{lcl}
 \theta_{k-1}^2t_k^2\norms{x_k - x_{k-1}}^2 - \theta_k^2t_{k+1}^2\norms{x_{k+1} - x_k}^2 & \geq & 2\gamma t_k^2(1-\nu_k)\iprods{Gy_{k-1}, x_{k+1} - x_k} \vspace{1ex}\\
 & \geq & -\frac{\gamma t_k^2}{t_{k+1}}\left[\norms{Gy_{k-1}}^2 + \norms{x_{k+1} - x_k}^2\right].
\end{array}
}
Following the same argument as in the proof of $\norms{x_k - y_k}^2$, we can show that $\lim_{k\to\infty}t_k^2\norms{x_{k+1} - x_k}^2 = 0$, and hence, $\norms{x_{k+1} - x_k}^2 = \SmallO{1/k^2}$, which proves the second part in the first line of  \eqref{eq:NAG_scheme00_convergence2}.
Since $\gamma^2\norms{Gy_k}^2 = \norms{x_{k+1} - y_k}^2 \leq 2\norms{x_{k+1} - x_k}^2 + 2\norms{y_k - x_k}^2$, we also obtain $\norms{Gy_k}^2 = \SmallO{1/k^2}$.

Since $\norms{Gx_k}^2 \leq 2\norms{Gx_k - Gy_k} + 2\norms{Gy_k}^2 \leq 2L^2\norms{x_k - y_k}^2 + 2\norms{Gy_k}^2$, we obtain the fourth line of \eqref{eq:NAG_scheme00_convergence2} from the previous lines.

Now, we prove the last line of \eqref{eq:NAG_scheme00_convergence2}.
Since $y_{k+1} - y_k  = x_{k+1} - x_k + \theta_{k}(x_{k+1} - x_k) - \theta_{k-1}(x_k - x_{k-1}) - \gamma\left(\nu_{k} Gy_k - \nu_{k-1} Gy_{k-1}\right)$, we can bound
\myeqn{
\hspace{-4ex}
\arraycolsep=0.2em
\begin{array}{lcl}
\norms{y_{k+1} - y_k}^2 & \leq & 16\norms{x_{k+1} - x_k}^2 + 4\norms{x_k - x_{k-1}}^2 + 4\gamma^2\norms{Gy_k - Gy_{k-1}}^2 \vspace{1ex}\\
&& + {~}  4\gamma^2 \norms{Gy_{k-1}}^2.
\end{array}
\hspace{-5ex}
}
Here, we have used the facts that $\theta_k, \theta_{k-1}, \nu_k, \nu_{k-1} \in (0, 1)$ and $(\nu_{k} - \nu_{k-1})^2 < 1$.
Combining this  inequality and the first and second lines of \eqref{eq:NAG_scheme00_convergence2}, we obtain the last line of \eqref{eq:NAG_scheme00_convergence2}.

\noindent\textbf{\color{red}The following proof is a correction for the published version on COAP.}
{\color{blue}Finally, to prove the convergence of $\sets{x_k}$ and $\sets{y_k}$, we note that $\Vc_k$ is nonnegative and nonincreasing, it converges. 
Moreover, from \eqref{eq:NAG_scheme00_Lyapunov_func1}, we also have $\norms{x_k - y^{\star}}^2 \leq \Vc_k$.
Hence, we conclude that $\sets{\norms{x_k - y^{\star}}}$ is bounded, i.e., there exists $M > 0$ such that $\norms{x_k - y^{\star}} \leq M$ for all $k\geq 0$.

Since $\mu = 1$, we can also express $\Vc_k$ from \eqref{eq:NAG_scheme00_Lyapunov_func1} as follows:
\begin{equation*}
\begin{array}{lcl}
\Vc_k &= & a_k\norms{Gy_{k-1}}^2 + b_k\iprods{Gy_{k-1}, x_k - y_k} + 2\norms{x_k - y^{\star}}^2 \vspace{1ex}\\
&& + {~} 2t_k\iprods{x_k - y^{\star}, y_k - x_k} + t_k^2\norms{y_k - x_k}^2
\end{array}
\end{equation*}
Since $\norms{x_k - y^{\star}} \leq M$ and $\lim_{k\to\infty}t_k\norms{y_k - x_k} = 0$, we have $\vert t_k\iprods{x_k - y^{\star}, y_k - x_k} \vert \leq t_k\norms{y_k - x_k}M \to 0$ as $k\to\infty$. 
We conclude that $\lim_{k\to\infty}t_k\iprods{x_k - y^{\star}, y_k - x_k} = 0$.
Similarly, since $b_k = 2\gamma t_k(t_k-1)$ and $\lim_{k\to\infty}t_k\norms{Gy_{k-1}} = \lim_{k\to\infty}t_k\norms{y_k - x_k} = 0$, we can easily show that $\lim_{k\to\infty}b_k\iprods{Gy_{k-1}, x_k - y_k}  = 0$.
Using these limits, $\lim_{k\to\infty}a_k\norms{Gy_{k-1}}^2 = \gamma^2\lim_{k\to\infty}t_k(t_k-1)\norms{Gy_{k-1}} = 0$, and the existence of $\lim_{k\to\infty}\Vc_k$ into the last expression, we conclude that $\lim_{k\to\infty}\norms{x_k - y^{\star}}$ exists.

Furthermore, since $G$ is $\frac{1}{L}$-co-coercive, it is $L$-Lipschitz continuous.
Therefore, any limit point $y^{\star}$ of $\sets{x_k}$ is in $\zer{G}$, and thus $\lim_{k\to\infty}\norms{x_k - y^{\star}} = 0$ by Opial's lemma.
We conclude that $\sets{x_k}$ is convergent to $y^{\star}$.
Since $\norms{x_k - y_k}\to 0$, combining this fact and $\lim_{k\to\infty}\norms{x_k - y^{\star}} = 0$, we also have $\lim_{k\to\infty}y_k = y^{\star}$.
}
\Eproof
\end{proof}

From the result of Theorem~\ref{th:NAG_scheme00_convergence2}, we can derive the convergence of Halpern's fixed-point iteration \eqref{eq:HP_scheme00}, but under different choice of parameters.

\begin{corollary}\label{co:HP_convergence2}
Let $\sets{y_k}$ be generated by Halpern's fixed-point iteration \eqref{eq:HP_scheme00} using $\beta_k := \frac{\omega+1}{k+2\omega +2}$ and $\eta_k := \gamma(1 - \beta_k)$ for a fixed $\gamma \in \left(0, \frac{2}{L}\right)$ and $\omega > 2$.
Then, the following statements hold:
\myeq{eq:convergence_of_HP_iter2}{
\arraycolsep=0.2em
\left\{\begin{array}{lclc|}
\sum_{k=0}^{\infty}(k + \omega + 1)\norms{Gy_{k}}^2 < +\infty \quad  & \text{and} \quad  & \norms{Gy_k}^2 = \SmallO{\frac{1}{k^2}}, \vspace{1ex}\\
\sum_{k=0}^{\infty}(k + \omega)\norms{y_{k+1} - y_k}^2 < +\infty \quad & \text{and} \quad & \norms{y_{k+1}  - y_k}^2 = \SmallO{\frac{1}{k^2}}.
\end{array}\right.
}
Consequently, $\set{y_k}$ converges to $y^{\star} \in \zer{G}$.
\end{corollary}

\begin{proof}
As proved in Theorem~\ref{th:HP2Nes_scheme},  \eqref{eq:HP_scheme00} is equivalent to \eqref{eq:NAG_scheme00} provided that  $\theta_k = \frac{\beta_k(1-\beta_{k-1})}{\beta_{k-1}}$, $\nu_k = \frac{\beta_k}{\beta_{k-1}}$, and $\gamma_k := \frac{\eta_k}{1 - \beta_k}$.
Using the choice of $\beta_k$, $\nu_k$, and $\gamma_k$ in Theorem~\ref{th:NAG_scheme00_convergence2}, we can show that $\theta_k = \frac{\beta_k(1-\beta_{k-1})}{\beta_{k-1}} = \frac{k+1}{k+2\omega+2}$ and $\nu_k =  \frac{\beta_k}{\beta_{k-1}}  = \frac{k + \omega + 2}{k+ 2\omega + 2}$.
These relations lead to $\beta_k =  \frac{\omega + 1}{k+2\omega +2}$.
Moreover, since $\gamma_k = \frac{\eta_k}{1-\beta_k} = \gamma \in \left(0, \frac{2}{L}\right)$, we have $\eta_k = \gamma(1-\beta_k)$.
Consequently, \eqref{eq:convergence_of_HP_iter2} follows from \eqref{eq:NAG_scheme00_convergence2}.
\Eproof
\end{proof}

If we set $\omega = 0$, then we obtain $\beta_k = \frac{1}{k+2}$ as in Theorem~\ref{th:HP_scheme00_convergence}.
In this case, we have to set $\mu = 0$ in $\Vc_k$ from \eqref{eq:NAG_scheme00_Lyapunov_func1}, and hence only obtain $\norms{Gy_k}^2 = \BigO{1/k^2}$ convergence rate.
Note that other choices of parameters in Theorem~\ref{th:NAG_scheme00_convergence2} are possible, e.g., by changing $\mu$ and $\omega$. 
Here, we have not tried to optimize the choice of these parameters.
As shown in \cite[Proposition 4.11]{Bauschke2011} that $T$ is a non-expansive mapping if and only if $G := \Id - T$ is $\frac{1}{2}$-co-coercive. 
Therefore, we can obtain new convergence results on the residual norm $\norms{y_k - Ty_k}$  from Corollary~\ref{co:HP_convergence2} for a Halpern's fixed-point iteration scheme to approximate a fixed-point $y^{\star}$ of $T$.

\beforesec
\section{Application to Monotone Inclusions}\label{sec:mono_inclusion}
\aftersec
In this section, we present three applications of Theorem~\ref{th:HP_scheme00_convergence} and Theorem~\ref{th:NAG_scheme00_convergence2} to proximal-point, forward-backward splitting, and three-operator splitting methods.

\beforesubsec
\subsection{Monotone inclusions and solution characterization}
\aftersubsec
We consider the following monotone inclusion and its special cases:
\myeq{eq:3o_MI}{
0 \in Ay^{\star} + By^{\star} + Cy^{\star},
}
where $A, B : \R^p \rightrightarrows 2^{\R^p}$ are multivalued and maximally monotone operators, and $C : \R^p \to\R^p$ is a $\frac{1}{L}$-co-coercive operator.
Let $Q := A + B + C$ and we assume that $\zer{Q} :=  Q^{-1}(0) = \sets{ y^{\star} \in \R^p : 0 \in Ay^{\star} + By^{\star} + Cy^{\star}}$ is nonempty.
We will consider the following cases in this paper.
\begin{itemize}
\itemsep=0.0em
\item\textbf{Case 1.} If $A=0$ and $B=0$, then by overloading $G = C$, \eqref{eq:3o_MI} reduces to the co-coercive equation \eqref{eq:MI} studied in Section~\ref{sec:main_part1}.

\item\textbf{Case 2.} If $B = 0$ and $C = 0$, then \eqref{eq:3o_MI} reduces to $0 \in Ay^{\star}$.
Then, we will investigate the convergence of an accelerated proximal-point algorithm and the interplay between Halpern's fixed-point iteration and Nesterov's accelerated interpretations in Subsection~\ref{subsec:prox_point_method}.

\item\textbf{Case 3.} If $C = 0$, then \eqref{eq:3o_MI} reduces to $0 \in Ay^{\star} + By^{\star}$, also covers monotone VIPs.
We will investigate the convergence of an accelerated forward-backward splitting scheme in Subsection~\ref{subssec:FB_method} using our results in Section~\ref{sec:main_part1}.

\item\textbf{Case 4.} Finally, we will also investigate the convergence of an accelerated three-operator splitting scheme for solving \eqref{eq:3o_MI}, and its special case: the accelerated Douglas-Rachford splitting scheme in Subsection~\ref{subsec:3o_method}.  
\end{itemize}

In order to characterize solutions of \eqref{eq:3o_MI}, we recall the following two operators.
The first operator is the forward-backward residual mapping associated with \textbf{Case 3} of \eqref{eq:3o_MI}  (i.e. $0 \in Ay^{\star} + By^{\star}$):
\myeq{eq:FB_residual}{
G_{\lambda Q}y := \tfrac{1}{\lambda}\left(y - J_{\lambda A}(y - \lambda By)\right), 
}
where $B$ is single-valued, $Q := A + B$, and $J_{\lambda A}$ is the resolvent of $\lambda A$ for any $\lambda > 0$.
The following result is proved similarly to  \cite[Proposition 26.1]{Bauschke2011} and \cite{tran2021halpern}.  

\begin{lemma}\label{leq:FBR_properties}
Let $A$ and $B$ in \eqref{eq:3o_MI} be maximally monotone, $B$ be single-valued, and $C = 0$.
Let $G_{\lambda Q}$ be defined by \eqref{eq:FB_residual}.
Then
\myeq{eq:DR:pro1}{
\iprods{G_{\lambda Q}x - G_{\lambda Q}y, x - y + \lambda (Bx - By)} \geq \lambda\norms{G_{\lambda Q}x - G_{\lambda Q}y}^2 + \iprods{Bx - By, x - y}.
}
Moreover, $G_{\lambda Q}y^{\star} = 0$ iff $y^{\star} \in \zer{A + B}$.
If, additionally, $B$ is $\frac{1}{L}$-co-coercive, then $G_{\lambda Q}$ is $\frac{\lambda(4 - \lambda L)}{4}$-co-coercive provided that  $0 < \lambda < \frac{4}{L}$.
\end{lemma}

Since $G_{\lambda Q}$ is co-coercive, Lemma \ref{leq:FBR_properties} shows that solving \eqref{eq:3o_MI} is equivalent to solving the co-coercive equation $G_{\lambda Q}y^{\star} = 0$ as a special case of  \eqref{eq:MI}.

The second operator is the residual mapping of a three-operator splitting scheme, which is defined as 
\myeq{eq:DR_residual}{
E_{\lambda Q}y  :=  \tfrac{1}{\lambda}(J_{\lambda B}y - J_{\lambda A}( 2J_{\lambda B}y - y - \lambda C\circ J_{\lambda B}y)),
}
where $J_{\lambda A} $ and $J_{\lambda B}$ are the resolvents of $\lambda A$ and $\lambda B$, respectively, and $\circ$ is a composition operator.
The following result is similar to the one in, e.g., \cite{Bauschke2011,Davis2015,he2015convergence} and we omit its proof here.

\begin{lemma}\label{le:DR:E_properties}
Let $A$ and $B$ in \eqref{eq:3o_MI} be maximally monotone, and $C$ be $\frac{1}{L}$-co-coercive. 
Let $E_{\lambda Q}$ be defined by \eqref{eq:DR_residual}. 
Then, $E_{\lambda Q}u^{\star} = 0$ iff $y^{\star} \in \zer{A+B+C}$, where $y^{\star} = J_{\lambda B}u^{\star}$.
Moreover, $E_{\lambda Q}$ satisfies the following property for all $u$ and $v$:
\myeq{eq:DR:E_pro1}{
\iprods{E_{\lambda Q}u - E_{\lambda Q}v, u - v} \geq \tfrac{\lambda(4 - L\lambda)}{4}\norms{E_{\lambda Q}u - E_{\lambda Q}v}^2.
}
If $B$ is single-valued and $C = 0$, then we have $E_{\lambda Q}u = G_{\lambda Q}y$, where $y = J_{\lambda B}u$ $($or equivalently, $u = y + \lambda By$$)$.
\end{lemma}

\beforesubsec
\subsection{Application to proximal-point method}\label{subsec:prox_point_method}
\aftersubsec
We consider \textbf{Case 2} where \eqref{eq:3o_MI} reduces to finding $y^{\star}\in\R^p$ such that $0 \in Ay^{\star}$.
Let $J_{\lambda A}y := (\Id + \lambda A)^{-1}y$ be the resolvent of $\lambda A$ for any $\lambda > 0$ and $G_{\lambda A}y = \frac{1}{\lambda}(\Id - J_{\lambda A})y = \frac{1}{\lambda}(y - J_{\lambda A}y)$ be the Yosida approximation of $A$ with index $\lambda > 0$.
Then, by  \cite[Corollary 23.11]{Bauschke2011}, $G_{\lambda A}$ is $\lambda$-co-coercive.
Moreover, $y^{\star}$ solves $0 \in Ay^{\star}$ if and only if $G_{\lambda A}y^{\star} = 0$.
Hence, solving $0 \in Ay^{\star}$ is equivalent to solving the $\lambda$-co-coercive equation $G_{\lambda A}y^{\star} = 0$.

In this case, the Halpern-type fixed-point scheme \eqref{eq:HP_scheme00} applying to $G_{\lambda A}y^{\star} = 0$, or equivalently, to solving $0 \in Ay^{\star}$, can be written as
\myeq{eq:HP_prox_point_scheme}{
y_{k+1} := \beta_ky_0 + (1 - \beta_k)y_k - \eta_k G_{\lambda A}y_k = \beta_ky_0 + \left(1 - \beta_k - \frac{\eta_k}{\lambda}\right)y_k + \frac{\eta_k}{\lambda}J_{\lambda A}y_k,
}
where $\beta_k$ and $\eta_k$ can be chosen either in Theorem~\ref{th:HP_scheme00_convergence} or Corollary~\ref{co:HP_convergence2} to guarantee convergence of \eqref{eq:HP_prox_point_scheme}.
If $\beta_k := \frac{1}{k+2}$ and $\eta_k := 2\lambda(1-\beta_k)$ as in Theorem \ref{th:HP_scheme00_convergence}, then  
\myeqn{
\arraycolsep=0.2em
\begin{array}{lcl}
y_{k+1} &:= & \beta_ky_0 + (1-\beta_k)y_k - 2(1-\beta_k)(y_k - J_{\lambda A}y_k) = \beta_ky_0 + (1-\beta_k)R_{\lambda A}y_k,
\end{array}
}
where $R_{\lambda A} := 2J_{\lambda A} - \Id$ is the reflected resolvent of $\lambda A$.
Moreover, under this choice of parameters, we have the following result from Theorem~\ref{th:HP_scheme00_convergence}:
\myeqn{
\norms{G_{\lambda A}y_k} \leq \frac{\norms{y_0 - y^{\star}}}{\lambda(k+1)}.
}
If we choose  $\beta_k := \frac{\omega + 1}{k + 2\omega + 2}$ and $\eta_k := \gamma(1 - \beta_k)$ as in  Corollary~\ref{co:HP_convergence2}, then \eqref{eq:HP_prox_point_scheme} becomes
\myeqn{
y_{k+1} := \frac{\omega+1}{k+2\omega+2} \cdot y_0 + \frac{k+\omega+1}{k+2\omega+2} \cdot \left[ \left(1 - \frac{\gamma}{\lambda} \right)y_k + \frac{\gamma}{\lambda} J_{\lambda A}y_k \right].
}
This expression can be viewed as  a new variant of Halpern's fixed-point iteration applied to the averaged mapping $\mathcal{T}_{\rho A}y = (1 - \rho)y + \rho J_{\lambda A}y$ with $\rho := \frac{\gamma}{\lambda}$ provided that $\gamma \in (0, \lambda]$.
In this case, we obtain a convergence result as in \eqref{eq:convergence_of_HP_iter2}.

Alternatively, if we apply \eqref{eq:NAG_scheme00} to solve $G_{\lambda A}y^{\star} = 0$, then we obtain a Nesterov's accelerated interpretation of \eqref{eq:HP_prox_point_scheme} as
\myeq{eq:NAG_scheme00_prox_point}{
\arraycolsep=0.2em
\left\{\begin{array}{lcl}
x_{k+1} &:= & y_k - \gamma_k G_{\lambda A}y_k = \left(1 - \rho_k \right)y_k +  \rho_k J_{\lambda A}y_k \quad \text{with} \quad \rho_k := \frac{\gamma_k}{\lambda}, \vspace{1ex}\\
y_{k+1} &:= & x_{k+1} + \theta_k(x_{k+1} - x_k) + \nu_k(y_k - x_{k+1}).
\end{array}\right.
}
This method was studied in \cite{mainge2021accelerated}.
Nevertheless, our analysis in Theorem~\ref{th:NAG_scheme00_convergence2} is simpler than that of \cite{mainge2021accelerated} when it applies to \eqref{eq:NAG_scheme00_prox_point}.
In particular, if we choose $\gamma_k := \lambda$, then the first line of \eqref{eq:NAG_scheme00_prox_point} reduces to $x_{k+1} = J_{\lambda A}y_k$.
The convergence rate guarantees of \eqref{eq:NAG_scheme00_prox_point} can be obtained as results of  Corollary~\ref{co:HP_scheme00_convergence} and Theorem~\ref{th:NAG_scheme00_convergence2}, respectively.

Finally, if we apply \eqref{eq:NAG_scheme00_2corr} to solve $G_{\lambda A}y^{\star} = 0$ and choose $\eta_k := \lambda\left(\frac{\beta_k}{\beta_{k-1}} + 1 - \beta_k \right)$ and $\gamma_k := \lambda$ such that $\nu_k = 0$ and $\kappa_k = \frac{\beta_k}{\beta_{k-2}}$, then \eqref{eq:NAG_scheme00_2corr} reduces to
\myeqn{
\arraycolsep=0.2em
\left\{\begin{array}{lcl}
x_{k+1} &:= & J_{\lambda A}y_k, \vspace{1ex}\\
y_{k+1} &:= & x_{k+1} + \theta_k(x_{k+1} - x_k) + \kappa_k(y_{k-1} - x_k).
\end{array}\right.
}
Clearly, if we choose $\beta_k := \frac{1}{k+2}$, then $\theta_k = \frac{k}{k+2}$ and $\kappa_k = \frac{k}{k+2}$.
This scheme reduces to the accelerated proximal-point algorithm in \cite{kim2021accelerated}.
In addition, we have $\eta_k = \frac{2\lambda(k+1)}{k+2}$ as in Theorem~\ref{th:HP_scheme00_convergence}.
Hence, the result of Corollary~\ref{co:HP_scheme00_convergence} is still applicable to this scheme to obtain a convergence rate guarantee $\norms{G_{\lambda Q}y_k} \leq \frac{\norms{y_0 - y^{\star}}}{\lambda(k+1)}$ as in \cite[Theorem 4.1]{kim2021accelerated}.

\beforesubsec
\subsection{Application to forward-backward splitting method}\label{subssec:FB_method}
\aftersubsec
Let us consider \textbf{Case 3} when \eqref{eq:3o_MI} reduces to finding $y^{\star}\in\R^p$ such that $0 \in  Ay^{\star} + By^{\star}$.
By Lemma~\ref{leq:FBR_properties}, $y^{\star} \in \zer{A+B}$ if and only if $G_{\lambda Q}y^{\star} = 0$, where $Q := A + B$ and $G_{\lambda Q}$ is defined by \eqref{eq:FB_residual}.
Moreover, $G_{\lambda Q}$ is $\frac{\lambda(4 - \lambda L)}{4}$-co-coercive, provided that $0 < \lambda < \frac{4}{L}$.

If we apply \eqref{eq:HP_scheme00} to solve $G_{\lambda Q}y^{\star} = 0$, then its iterate can be written as
\myeq{eq:HP_FB_scheme}{
y_{k+1} := \beta_ky_0 + (1-\beta_k) \big[  (1 - \rho) y_k + \rho J_{\lambda A}(y_k - \lambda By_k) \big],
}
where we have set $\rho :=  \frac{4-\lambda L}{2}$.
In particular, if we choose $\lambda := \frac{2}{L}$, then $\rho = 1$ and \eqref{eq:HP_FB_scheme} reduces to $y_{k+1} := \beta_ky_0 + (1-\beta_k)J_{\lambda A}(y_k - \lambda By_k)$, which can be viewed as Halpern's fixed-point iteration applied to approximate a fixed-point of $J_{\lambda A}(y_k - \lambda By_k)$.

Depending on the choice of $\beta_k$ and $\rho$ as in Theorem~\ref{th:HP_scheme00_convergence} or Corollary~\ref{co:HP_convergence2}, we obtain 
\myeqn{
\norms{G_{\lambda Q}y_k} \leq \frac{4\norms{y_0 - y^{\star}}}{\lambda(4 - \lambda L)(k+1)}, \quad \text{or} \quad \norms{G_{\lambda Q}y_k} = \SmallO{1/k}, 
}
respectively, provided that $0 < \lambda < \frac{4}{L}$. 

Now, we consider Nesterov's accelerated variant of \eqref{eq:HP_FB_scheme} by applying \eqref{eq:NAG_scheme00} to $G_{\lambda Q}y^{\star} = 0$ to obtain the following one:
\myeq{eq:NES_FB_scheme}{
\arraycolsep=0.3em
\left\{\begin{array}{lcl}
x_{k+1} &:= & (1 - \rho_k) y_k + \rho_k J_{\lambda A}(y_k - \lambda By_k), \vspace{1ex}\\
y_{k+1} &:= & x_{k+1} + \theta_k(x_{k+1} - x_k) + \nu_k(y_k - x_{k+1}), 
\end{array}\right.
}
where $\rho_k := \frac{\gamma_k}{\lambda}$, $\theta_k :=   \frac{\beta_k(1-\beta_{k-1})}{\beta_{k-1}}$, and $\nu_k := \frac{\beta_k}{\beta_{k-1}}$.
This scheme is similar to the one studied in \cite{mainge2021fast}.
Again, the convergence of \eqref{eq:NES_FB_scheme} can be guaranteed by either Corollary~\ref{co:HP_scheme00_convergence} or Theorem~\ref{th:NAG_scheme00_convergence2} depending on the choice of $\gamma_k$, $\theta_k$, and $\nu_k$.
However, we omit the details here.

\beforesubsec
\subsection{Application to three-operator splitting method}\label{subsec:3o_method}
\aftersubsec
Finally, we consider the general case, \textbf{Case 4}.
As stated in Lemma~\ref{le:DR:E_properties}, $y^{\star} \in \zer{A+B+C}$ if and only if $E_{\lambda Q}y^{\star} = 0$, where $Q := A + B + C$ and $E_{\lambda Q}$ is defined by \eqref{eq:DR_residual}.
Let us apply \eqref{eq:HP_scheme00}  to $E_{\lambda Q}y^{\star} = 0$ and arrive at the following scheme:
\myeqn{
\arraycolsep=0.3em
\begin{array}{lcl}
y_{k+1} & := & \beta_ky_0 + (1-\beta_k) y_k - \frac{\eta_k}{\lambda}(J_{\lambda B}y_k - J_{\lambda A}( 2J_{\lambda B}y_k - y_k - \lambda C \circ J_{\lambda B}y_k)).
\end{array}
}
Unfolding this scheme by using intermediate variables $z_k$ and $w_k$, we get
\myeq{eq:DR_HT_scheme}{
\arraycolsep=0.3em
\left\{\begin{array}{lcl}
z_k &:= & J_{\lambda B}y_k, \vspace{1ex}\\
w_k &:= & J_{\lambda A}(2z_k - y_k - \lambda Cz_k), \vspace{1ex}\\
y_{k+1} & := & \beta_ky_0 + (1-\beta_k)y_k - \frac{\eta_k}{\lambda} (z_k - w_k).
\end{array}\right.
}
This is called a Halpern-type three-operator splitting scheme for solving \eqref{eq:3o_MI}.
If $C = 0$, then it reduces to a Halpern-type  Douglas-Rachford splitting scheme for solving \textbf{Case 3} of \eqref{eq:3o_MI} derived from \eqref{eq:HP_scheme00}.
The latter case was proposed in \cite{tran2021halpern} with a direct convergence proof for both dynamic and constant stepsizes, but the convergence is given on $G_{\lambda Q}$ instead of $E_{\lambda Q}$.
Note that the convergence results of Theorem~\ref{th:HP_scheme00_convergence} and Corollary~\ref{co:HP_convergence2} can be applied to \eqref{eq:DR_HT_scheme} to obtain convergence rates on $\norms{E_{\lambda Q}y_k}$.
Such rates can be transformed into the ones on $\norms{G_{\lambda Q}z_k}$ when $C = 0$ and $B$ is single-valued.

Next, we can also derive Nesterov's accelerated variant of \eqref{eq:DR_HT_scheme} by applying \eqref{eq:NAG_scheme00} to solve $E_{\lambda Q}y^{\star} = 0$.
In this case, \eqref{eq:NAG_scheme00} becomes 
\myeq{eq:DR_NAG_scheme}{
\arraycolsep=0.3em
\left\{\begin{array}{lcl}
z_k &:= & J_{\lambda B}y_k, \vspace{1ex}\\
w_k &:= & J_{\lambda A}(2z_k - y_k - \lambda Cz_k), \vspace{1ex}\\
x_{k+1} &:= & y_k + \frac{1}{\lambda}(w_k - z_k), \vspace{1ex}\\
y_{k+1} &:= & x_{k+1} + \theta_k(x_{k+1} - x_k) + \nu_k(y_k - x_{k+1}).
\end{array}\right.
}
Here, the parameters $\theta_k$ and $\nu_k$ can be chosen as in either Corollary~\ref{co:HP_scheme00_convergence} or Theorem~\ref{th:NAG_scheme00_convergence2}.
This scheme essentially has the same per-iteration complexity as the standard three-operator splitting scheme in the literature, including \cite{Davis2015}.
However, its convergence rate is much faster than the standard one by applying either Corollary~\ref{co:HP_scheme00_convergence} or Theorem~\ref{th:NAG_scheme00_convergence2}. 
If $C = 0$, then \eqref{eq:DR_NAG_scheme} reduces to an accelerated Douglas-Rachford splitting scheme, where its fast convergence rate can be obtained as a special case of either Corollary~\ref{co:HP_scheme00_convergence} or Theorem~\ref{th:NAG_scheme00_convergence2}. 

\beforesec
\section{Extra-Anchored Gradient Method and Its Variants}\label{sec:EAG_sec}
\aftersec
\paragraph{Motivation.}
While the gradient of a convex and $L$-smooth function is co-coercive, monotone and Lipschitz continuous operators are not co-coercive in general.
As a simple example, one can take $Gx = (Av, -A^{\top}u)$ as the gradient of  the saddle objective function in a bilinear game, where $A$ in $\R^{m\times n}$ is given and $x = (u, v)$.
In order to solve \eqref{eq:MI} when $G$ is only monotone and $L$-Lipschitz continuous, the extragradient method (EG) appears to be one of the most suitable candidates \cite{Korpelevic1976}. 
This method has recently been extended to weak Minty VIP, i.e. $\iprods{Gy, y - y^{\star}} \geq -\rho\norms{Gy}^2$ for all $y\in\R^p$ and $y^{\star} \in \zer{G}$ in our context, see, e.g. \cite{diakonikolas2021efficient}.
In \cite{yoon2021accelerated}, Yoon and Ryu applied Halpern's fixed-point iteration to EG and obtained a new algorithm called extra-anchored gradient method (EAG).
This algorithm achieves optimal convergence rate on $\norms{Gy_k}$.
Recently,  \cite{lee2021fast} extended EAG to a co-monotone setting of \eqref{eq:MI} and still achieved the same rate $\norms{Gy_k} = \BigO{1/k}$ as in \cite{yoon2021accelerated}.
This is perhaps surprising since $G$ is nonmonotone. 
An extension to Popov's scheme (also called past-extra-gradient, or reflected forward methods) can be found in \cite{tran2021halpern}.
Our goal in this section is to derive a corresponding Nesterov's accelerated interpretation of these schemes and possibly provide an alternative convergence rate analysis for the existing results in  \cite{lee2021fast,tran2021halpern,yoon2021accelerated}.

\beforesubsec
\subsection{The extra-anchored gradient method and its convergence}
\aftersubsec
The extra-anchored gradient method (EAG) was proposed in \cite{yoon2021accelerated} to solve \eqref{eq:MI} under the monotonicity and $L$-Lipschitz continuity of $G$, which can be written as
\myeq{eq:EAG_scheme30}{
\arraycolsep=0.3em
\left\{\begin{array}{lcl}
z_{k+1} &:= & \beta_ky_0 + (1-\beta_k)y_k - \eta_kGy_k, \vspace{1ex}\\
y_{k+1} &:= & \beta_ky_0 + (1-\beta_k)y_k - \hat{\eta}_kGz_{k+1}, 
\end{array}\right.
}
where $\beta_k \in (0, 1)$, $y_0$ is an initial point, and $\eta_k$ and $\hat{\eta}_k$ are two given step-sizes.
Here, we use two different step-sizes $\eta_k$ and $\hat{\eta}_k$ compared to the original EAG in \cite{yoon2021accelerated} by adopting the idea of EG+ from \cite{diakonikolas2021efficient}, see also \cite{lee2021fast}.

As proven in \cite{yoon2021accelerated}, if we update $\beta_k := \frac{1}{k+2}$, $ \eta_{k+1} := \big(1 - \frac{L^2\eta_k^2}{(1 - L^2\eta_k^2)(k+1)(k+3)} \big)\eta_k$, and $\hat{\eta}_{k+1}  := \eta_{k+1}$ with some $0 < \eta_0 < \frac{1}{L}$, then we obtain
\myeq{eq:p7:EAG_convergence}{
\hspace{-1ex}
\norms{Gy_k}^2 \leq  \frac{C_{*}\norms{y_0 - y^{\star}}^2}{(k+1)(k+2)}, \quad \text{where} \ \eta_{*} = \lim_{k\to\infty}\eta_k > 0 \ \text{and} \ C_{*} := \tfrac{4(1 + \eta_0\eta_{*}L^2)}{\eta_{*}^2}.
\hspace{-1ex}
}
Alternatively, one can also fix the step-size $\hat{\eta}_k = \eta_k = \eta \in \left(0, \frac{1}{8L}\right]$, and the following convergence guarantee is established:
\myeq{eq:p7:EAG_convergence_const}{
\norms{Gy_k}^2 \leq \frac{C_{*}\norms{y_0 - y^{\star}}^2}{(k+1)^2}, \quad\text{where} \quad C_{*} :=  \frac{4(1 + \eta L + \eta^2L^2)}{\eta^2(1 + \eta L)}.
}
In particular, if $\eta := \frac{1}{8L}$, then $C_{*} = 260$.
Both \eqref{eq:p7:EAG_convergence} and \eqref{eq:p7:EAG_convergence_const} were proven in \cite{yoon2021accelerated}.

\beforesubsec
\subsection{Nesterov's accelerated interpretation of EAG}\label{subsec:EAG2NAG_scheme}
\aftersubsec
Now, let us derive  Nesterov's accelerated interpretation of \eqref{eq:EAG_scheme30} by proving the following result.

\begin{theorem}\label{th:HP2Nes_scheme_EAG}
Let $\sets{(x_k, y_k, z_k)}$ be generated by the following scheme:
\myeq{eq:NAG_scheme10}{
\arraycolsep=0.3em
\left\{\begin{array}{lcl}
x_{k+1} & := & y_k - \gamma_k Gy_k, \vspace{1ex}\\
z_{k+1} & := & x_{k+1}  +  \theta_k(x_{k+1} - x_k) +  \nu_k(z_k - x_{k+1}), \vspace{1ex}\\
y_{k+1} & := & z_{k+1} - \hat{\eta}_kGz_{k+1} + \eta_kGy_k, 
\end{array}\right.
}
starting from $z_0 = x_0 := y_0$, where $\gamma_k := \frac{\eta_k}{1-\beta_k}$, $\theta_k := \frac{\beta_k(1-\beta_{k-1})}{\beta_{k-1}}$, and $\nu_k := \frac{\beta_k}{\beta_{k-1}}$ for $\beta_k$, $\eta_k$, and $\hat{\eta}_k$ given in \eqref{eq:EAG_scheme30}.
Then, $\sets{(y_k, z_k)}$ is identical to the one generated by the EAG scheme \eqref{eq:EAG_scheme30}.
\end{theorem}

\revise{Clearly, \eqref{eq:NAG_scheme10} is new compared to any Nesterov's accelerated scheme in the literature.}
To see a relation to existing methods, we rewrite \eqref{eq:NAG_scheme10} equivalently to
\myeq{eq:NAG_scheme10b}{
\arraycolsep=0.3em
\left\{\begin{array}{lcl}
x_{k+1} & := & z_k - \hat{\eta}_{k-1}\hat{G}z_k, \vspace{1ex}\\
z_{k+1} & := & x_{k+1}  +  \theta_k(x_{k+1} - x_k) +  \nu_k(z_k - x_{k+1}), 
\end{array}\right.}
where $\hat{G}z_k :=  Gz_k + \frac{\gamma_k}{\hat{\eta}_{k-1}} Gy_k  - \frac{\eta_{k-1}}{\hat{\eta}_{k-1}}Gy_{k-1}$ and $y_{k+1} :=  z_{k+1} - \hat{\eta}_kGz_{k+1} + \eta_kGy_k$.
Clearly, the first two lines of \eqref{eq:NAG_scheme10b} are similar to \eqref{eq:NAG_scheme00}, but using an approximate operator $\hat{G}z_k$ instead of the exact evaluation $Gz_k$ as in \eqref{eq:NAG_scheme00}.
Therefore, \eqref{eq:NAG_scheme10b} can be viewed as an inexact variant of the Nesterov's accelerated method \eqref{eq:NAG_scheme00}.

\begin{proof}[\textbf{Proof of Theorem~\ref{th:HP2Nes_scheme_EAG}}]
We only prove that \eqref{eq:EAG_scheme30} leads to  \eqref{eq:NAG_scheme10}.
The opposite direction from \eqref{eq:NAG_scheme10}  to  \eqref{eq:EAG_scheme30} is obtained by reverting back the derivations below.  

\revise{\noindent Firstly, multiplying the first line of \eqref{eq:EAG_scheme30} by $\beta_{k-1}$, we have 
\begin{equation*}
\beta_{k-1}z_{k+1} = \beta_k\beta_{k-1}y_0 + \beta_{k-1}(1 - \beta_k)y_k - \beta_{k-1}\eta_kGy_k.
\end{equation*}
Shifting the index from $k$ to $k-1$ of the first line of \eqref{eq:EAG_scheme30}, and then multiplying the result by $-\beta_k$, we get 
\begin{equation*}
-\beta_kz_k = -\beta_k\beta_{k-1}y_0 - \beta_k(1-\beta_{k-1})y_{k-1} + \beta_k\eta_{k-1}Gy_{k-1}.
\end{equation*}
Summing up both expressions, we arrive at 
\begin{equation*}
\beta_{k-1}z_{k+1} - \beta_{k}z_k = \beta_{k-1}(1-\beta_k)y_k - \beta_{k-1}\eta_kGy_k - \beta_k(1-\beta_{k-1})y_{k-1} + \beta_k\eta_{k-1}Gy_{k-1}.
\end{equation*}
This expression leads to 
\begin{equation*}
z_{k+1} = \frac{\beta_k}{\beta_{k-1}}z_k + (1-\beta_{k})y_{k} - \eta_{k}Gy_{k} - \frac{\beta_{k}(1-\beta_{k-1})}{\beta_{k-1}}y_{k-1} + \frac{\beta_{k}\eta_{k-1}}{\beta_{k-1}}Gy_{k-1}.
\end{equation*}
}
Next, subtracting the first line from the second one of \eqref{eq:EAG_scheme30}, we have $y_{k+1} - z_{k+1} = - \hat{\eta}_kGz_{k+1} +  \eta_k Gy_k$, leading to $y_{k+1} = z_{k+1} - \hat{\eta}_kGz_{k+1} + \eta_kGy_k$.
Combining this expression and the last line above, we obtain
\myeq{eq:AEG_iter2}{
\hspace{-2ex}
\arraycolsep=0.3em
\left\{\begin{array}{lcl}
z_{k+1} & = & \frac{\beta_k}{\beta_{k-1}}z_k + (1-\beta_{k})y_{k} - \eta_{k}Gy_{k} - \frac{\beta_{k}(1-\beta_{k-1})}{\beta_{k-1}}y_{k-1} + \frac{\beta_{k}\eta_{k-1}}{\beta_{k-1}}Gy_{k-1}, \vspace{1ex}\\
y_{k+1} & = & z_{k+1} - \hat{\eta}_kGz_{k+1} + \eta_k Gy_k.
\end{array}\right.
\hspace{-2ex}
}
Let us introduce $x_{k+1} := y_k - \gamma_k Gy_k$.
Then, we have $Gy_k = \frac{1}{\gamma_k}(y_k - x_{k+1})$.
Substituting these expressions into the first line of \eqref{eq:AEG_iter2}, we get
\myeqn{ 
\hspace{-2ex}
\arraycolsep=0.1em
\begin{array}{lcl}
z_{k+1} & = &  \frac{\beta_{k}}{\beta_{k-1}}z_{k} + (1-\beta_{k})y_{k} -  \frac{\eta_k}{\gamma_k}(y_k - x_{k+1}) - \frac{\beta_{k}(1-\beta_{k-1})}{\beta_{k-1}}y_{k-1} + \frac{\beta_{k}\eta_{k-1}}{\beta_{k-1}\gamma_{k-1}}(y_{k-1} - x_k) \vspace{1ex}\\
&= & \frac{\beta_k}{\beta_{k-1}}z_{k}  + \left(1 - \frac{\beta_k}{\beta_{k-1}}\right)x_{k+1} +  \frac{\beta_k(1-\beta_{k-1})}{\beta_{k-1}}(x_{k+1} - x_k)  \vspace{1ex}\\
&& + {~}  \left(1 - \beta_k - \frac{\eta_k}{\gamma_k}\right)(y_k - x_{k+1}) -  \frac{\beta_k}{\beta_{k-1}}\big( 1-\beta_{k-1} - \frac{\eta_{k-1}}{\gamma_{k-1}} \big)(y_{k-1} - x_k).
\end{array}
\hspace{-2ex}
}
If we choose $\gamma_k$ such that $1 - \beta_k - \frac{\eta_k}{\gamma_k} = 0$ (or equivalently, $\gamma_k = \frac{\eta_k}{1 - \beta_k}$), then we have 
\myeqn{
\arraycolsep=0.3em
\begin{array}{lcl}
z_{k+1} & = & \frac{\beta_k}{\beta_{k-1}}z_{k}  + \left(1 - \frac{\beta_k}{\beta_{k-1}}\right)x_{k+1} +  \frac{\beta_k(1-\beta_{k-1})}{\beta_{k-1}}(x_{k+1} - x_k).
\end{array}
}
Finally, putting the above derivations together, we eventually get \eqref{eq:NAG_scheme10}.
\Eproof
\end{proof}

In order to analyze the convergence of \eqref{eq:NAG_scheme10}, following the same approach as in Section \ref{sec:main_part1}, we consider the following Lyapunov function:
\myeq{eq:NAG_scheme10_Lyapunov}{
\Qc_k := a_k\norms{Gy_{k-1}}^2 + b_k\iprods{Gy_{k-1}, x_k - z_k} + \norms{x_k + t_k(z_k - x_k) - y^{\star} }^2 + \mu\norms{x_k - y^{\star}}^2,
}
where $a_k > 0$, $b_k > 0$, $t_k > 0$, and $\mu \geq 0$ are given, determined later.

\revise{Recall that the Lyapunov function used in \cite{yoon2021accelerated} is $\Lc_k := p_k\norms{Gy_k} + q_k\iprods{Gy_k, y_k - y_0}$, where $p_k$ and $q_k$ are given.
Similar to Lemma~\ref{le:from_L_to_V}, we can show that if $\mu = 0$,  $a_k := \frac{4p_k(p_k + \gamma_kq_k)}{q_k^2}$, and $b_k := \frac{4p_k}{\beta_kq_k}$, then $\Lc_k = \frac{q_k^2}{4p_k}\left[ \Qc_{k+1} - \norms{y_0 - y^{\star}}^2\right]$.
This relation allows one to adopt the analysis in \cite{yoon2021accelerated} to prove convergence of \eqref{eq:NAG_scheme10}.
However, if  $\frac{q_k^2}{4p_k}$ is not a constant, then $\Qc_{k+1}$ remains different from $\Lc_k$.}

The following lemma proves a monotone property of $\Qc_k$, which plays a key role to establish convergence of \eqref{eq:NAG_scheme10}.

\begin{lemma}\label{le:NAG_scheme10_descent}
\revise{
Suppose that $G$ in \eqref{eq:MI} is monotone and $L$-Lipschitz continuous.
Let $\sets{(x_k, y_k, z_k)}$ be generated by \eqref{eq:NAG_scheme10} and $\Qc_k$ be defined by \eqref{eq:NAG_scheme10_Lyapunov} with $\mu := 0$.
For $t_k > 1$, assume that the parameters $\gamma_k$, $\hat{\eta}_k$, $\eta_k$,  $\theta_k$, $\nu_k$, $a_k$, and $b_k$ are updated by
\myeq{eq:NAG_scheme10_param}{
\hspace{-2ex}
\arraycolsep=0.2em
\begin{array}{ll}
& \gamma_k := \gamma \in \left(0,  \frac{1}{L}\right], \quad  \hat{\eta}_k :=  \gamma, \quad \eta_k := \frac{\gamma(t_{k+1} - 1)}{t_{k+1}}, \quad \theta_k := \frac{t_k-1}{t_{k+1}}, \quad \nu_k := \frac{t_k}{t_{k+1}}, \vspace{1ex}\\
&a_k := \frac{\gamma b_k(t_{k-1} + 2)}{2t_k}, \quad\text{and}\quad b_{k+1} := \frac{b_kt_{k+1}}{t_k-1}.
\end{array} 
\hspace{-2ex}
}
Then, for $k\geq 0$, the following estimates hold:
\myeq{eq:NAG_scheme10_descent}{
\arraycolsep=0.2em
\begin{array}{lcl}
\Qc_k - \Qc_{k+1} & \geq & \frac{b_k(1 - L^2\gamma^2)t_k}{2L^2\gamma(t_k-1)} \norms{Gy_k - Gz_k}^2 + \frac{\gamma b_k(t_{k-1} - t_k + 1)}{2t_k}\norms{Gy_{k-1}}^2, \vspace{1ex}\\
\Qc_{k+1} & \geq & \frac{\gamma b_k}{2t_k} \left( t_{k-1} - \frac{b_k}{2\gamma t_k} \right)\norms{Gy_{k-1}}^2.
\end{array}
}
}
\end{lemma}
%
\begin{proof}
Similar to the proof of \eqref{eq:Vk_Vk1_proof5}, using \eqref{eq:NAG_scheme10} and \eqref{eq:NAG_scheme10_Lyapunov}, we can prove that
\myeq{eq:NAG_scheme10_proof4}{
\arraycolsep=0.3em
\begin{array}{lcl}
\Qc_k - \Qc_{k+1} &= & a_k\norms{Gy_{k-1}}^2 - a_{k+1}\norms{Gy_k}^2 - b_k\iprods{Gy_{k-1}, x_{k+1} - x_k} \vspace{1ex}\\
&& + {~} b_{k+1}\theta_k\iprods{Gy_k, x_{k+1} - x_k}  - b_k\iprods{Gy_{k-1}, z_k - x_{k+1}} \vspace{1ex}\\
&& + {~} b_{k+1}\nu_k\iprods{Gy_k, z_k - x_{k+1}}  + (t_k^2 - \nu_k^2t_{k+1}^2)\norms{z_k - x_{k+1}}^2 \vspace{1ex}\\
&& + {~} \left[ (t_k-1)^2 - t_{k+1}^2\theta_k^2 + \mu \right]\norms{x_{k+1} - x_k}^2 \vspace{1ex}\\
&& + {~} 2(t_k - 1 - t_{k+1}\theta_k - \mu)\iprods{x_{k+1} - x_k, x_{k+1} - y^{\star}} \vspace{1ex}\\
&& + {~} 2(t_k - t_{k+1}\nu_k)\iprods{z_k - x_{k+1}, x_{k+1} - y^{\star}} \vspace{1ex}\\
&& + {~} 2\left[ t_k(t_k-1) - \nu_k\theta_kt_{k+1}^2 \right]\iprods{z_k - x_{k+1}, x_{k+1} - x_k}.
\end{array} 
}
Now, let us choose the parameters $t_k$, $b_k$, $\theta_k$, and $\nu_k$ such that 
\myeq{eq:NAG_scheme10_proof5}{
\hspace{-1ex}
\begin{array}{ll}
&t_k - t_{k+1}\nu_k = 0, \ t_k(t_k-1) - \nu_k\theta_kt_{k+1}^2 = 0, \vspace{1ex}\\
&  t_k - 1 - t_{k+1}\theta_k - \mu = 0, \ \text{and} \ b_k = b_{k+1}\theta_k.
\end{array}
\hspace{-2ex}
}
The third condition leads to $\theta_k := \frac{t_k - 1 - \mu}{t_{k+1}}$, while the first one gives us $\nu_k := \frac{t_k}{t_{k+1}}$.
The second condition becomes $t_k(t_k - 1) = t_k(t_{k-1} - 1 - \mu)$, which is satisfied if $\mu = 0$.
The last condition holds if  $b_{k+1} := \frac{b_k}{\theta_k}$.
These updates are exactly \eqref{eq:NAG_scheme10_param}.

Noticing  from \eqref{eq:NAG_scheme10} that $z_k - x_{k+1} = \gamma_kGy_k +  \hat{\eta}_{k-1}Gz_k - \eta_{k-1}Gy_{k-1}$.
Using this relation and our choice $b_{k+1} = \frac{b_k}{\theta_k}$, we can easily show that 
\myeqn{
\arraycolsep=0.3em
\begin{array}{lcl}
\Tc_{[2]} &:= & b_{k+1}\nu_k\iprods{Gy_k, z_k - x_{k+1}}  - b_k\iprods{Gy_{k-1}, z_k - x_{k+1}} \vspace{1ex}\\
&= & \frac{b_k\gamma_k\nu_k}{\theta_k}\norms{Gy_k}^2 - b_k\left(\gamma_k + \frac{\eta_{k-1}\nu_k}{\theta_k} \right)\iprods{Gy_k, Gy_{k-1}} + \frac{b_k\hat{\eta}_{k-1}\nu_k}{\theta_k}\iprods{Gy_k, Gz_k} \vspace{1ex}\\
&& - {~} b_k\hat{\eta}_{k-1}\iprods{Gy_{k-1}, Gz_k} + b_k\eta_{k-1}\norms{Gy_{k-1}}^2.
\end{array} 
}
Utilizing the monotonicity  of $G$, $b_k = b_{k+1}\theta_k$, and $x_{k+1} = y_k - \gamma_kGy_k$, we have
\myeqn{
\arraycolsep=0.1em
\begin{array}{lcl}
\Tc_{[3]} &:= & b_{k+1}\theta_k\iprods{Gy_k, x_{k+1} - x_k}  - b_k\iprods{Gy_{k-1}, x_{k+1} - x_k} \vspace{0.5ex}\\
&\geq & - b_k\gamma_k\norms{Gy_k}^2 + b_k(\gamma_k + \gamma_{k-1})\iprods{Gy_k, Gy_{k-1}} - b_k\gamma_{k-1}\norms{Gy_{k-1}}^2.
\end{array} 
}
Substituting \eqref{eq:NAG_scheme10_proof5}, $\Tc_{[2]}$, and $\Tc_{[3]}$ into \eqref{eq:NAG_scheme10_proof4} and using $\frac{\nu_k}{\theta_k} = \frac{t_k}{t_k-1-\mu} = \frac{t_k}{t_k-1}$, we can further lower bound
\myeq{eq:NAG_scheme10_proof6}{
\arraycolsep=0.3em
\begin{array}{lcl}
\Qc_k - \Qc_{k+1} &\geq & \left[ a_k + b_k( \eta_{k-1} - \gamma_{k-1}) \right] \norms{Gy_{k-1}}^2 + \left( \frac{b_k\gamma_k}{t_k-1}  - a_{k+1} \right) \norms{Gy_k}^2 \vspace{0.5ex}\\
&& - {~}  b_k\left(\frac{\eta_{k-1}t_k}{t_k-1} - \gamma_{k-1} \right)\iprods{Gy_k, Gy_{k-1}} \vspace{1ex}\\
&& + {~} \frac{b_k \hat{\eta}_{k-1}}{t_k-1} \iprods{Gy_k, Gz_k} +  b_k\hat{\eta}_{k-1}\iprods{Gy_k - Gy_{k-1}, Gz_k}.  
\end{array} 
}
Now, using the $L$-Lipschitz continuity of $G$, we have $\norms{Gz_k - Gy_k}^2 \leq L^2\norms{z_k - y_k}^2 = L^2\norms{\hat{\eta}_{k-1}Gz_k - \eta_{k-1}Gy_{k-1}}^2$, which leads to 
\myeqn{
\arraycolsep=0.1em
\begin{array}{lcl}
\norms{Gy_k}^2 & + &  (1 - L^2\hat{\eta}_{k-1}^2)\norms{Gz_k}^2 - 2\left(1 - L^2\eta_{k-1}\hat{\eta}_{k-1}\right) \iprods{Gy_k, Gz_k} \vspace{1ex}\\
&& - {~} 2L^2\eta_{k-1}\hat{\eta}_{k-1}\iprods{Gz_k, Gy_k - Gy_{k-1}} - L^2\eta_{k-1}^2\norms{Gy_{k-1}}^2 \leq 0.
\end{array} 
}
Multiplying this inequality by $\frac{b_k}{2L^2\eta_{k-1}}$ and adding the result to \eqref{eq:NAG_scheme10_proof6}, we get
\myeq{eq:NAG_scheme10_proof7}{
\hspace{-0ex}
\arraycolsep=0.2em
\begin{array}{lcl}
\Qc_k - \Qc_{k+1} &\geq & \frac{b_k(1 - L^2\hat{\eta}_{k-1}^2)}{2L^2\eta_{k-1}}\norms{Gz_k}^2  + \left( \frac{b_k\gamma_k}{t_k-1} + \frac{b_k}{2L^2\eta_{k-1}} - a_{k+1} \right) \norms{Gy_k}^2  \vspace{0.5ex}\\
&& - {~} b_k\left( \frac{1}{L^2\eta_{k-1}}  -  \frac{\hat{\eta}_{k-1}t_k }{t_k-1}   \right) \iprods{Gy_k, Gz_k} \vspace{1ex}\\
&& + {~} \left[ a_k + \frac{b_k}{2}( \eta_{k-1} - 2\gamma_{k-1})  \right] \norms{Gy_{k-1}}^2  \vspace{1ex}\\
&& - {~} b_k\left(\frac{\eta_{k-1}t_k}{t_k-1 } - \gamma_{k-1} \right)\iprods{Gy_k, Gy_{k-1}}.
\end{array} 
\hspace{-0ex}
}
\revise{Let us choose $\gamma_{k-1}$, $\hat{\eta}_{k-1}$, and $\eta_{k-1}$ as in \eqref{eq:NAG_scheme10_param}, i.e.:
\myeq{eq:NAG_scheme10_param2}{
\hspace{-2.0ex}
\arraycolsep=0.1em
\begin{array}{lcl}
\gamma_{k-1} = \hat{\eta}_{k-1} := \gamma \in \left(0, \frac{1}{L}\right], \quad \text{and} \quad \eta_{k-1} := \frac{\gamma(t_k-1)}{t_k}.
\end{array}
\hspace{-2ex}
}
From \eqref{eq:NAG_scheme10_proof5}, we have $b_k = b_{k+1}\theta_k = \frac{b_{k+1}(t_k - 1)}{t_{k+1}}$, leading to $b_{k+1} = \frac{b_kt_{k+1}}{t_k-1}$ as in \eqref{eq:NAG_scheme10_param}.
If we choose $a_{k+1} :=  \frac{\gamma b_k(t_k+2)}{2(t_k-1)}$, then since $b_{k-1} = \frac{b_k(t_{k-1}-1)}{t_k}$, we get $a_k = \frac{\gamma b_k(t_{k-1}+2)}{2t_k}$ as given in \eqref{eq:NAG_scheme10_param}.
Next, utilizing \eqref{eq:NAG_scheme10_param}, we can show that
\myeqn{
\arraycolsep=0.3em
\left\{\begin{array}{lcl}
\frac{b_k(1 - L^2\hat{\eta}_{k-1}^2)}{2L^2\eta_{k-1}} &= & \frac{b_k(1 - L^2\gamma^2)t_k}{2L^2\gamma(t_k-1)}, \vspace{0.5ex}\\
b_k\left( \frac{1}{L^2\eta_{k-1}} -  \frac{\hat{\eta}_{k-1} t_k }{t_k - 1} \right) &= & \frac{b_k(1 - L^2\gamma^2)t_k}{L^2\gamma(t_k-1)},  \vspace{1ex}\\
\frac{b_k\gamma_k}{t_k-1}  + \frac{b_k}{2L^2\eta_{k-1}} - a_{k+1} &= &  \frac{b_k(1 - L^2\gamma^2)t_k}{2L^2\gamma(t_k-1)},  \vspace{1ex}\\
a_k + \frac{b_k}{2}( \eta_{k-1} - 2\gamma_{k-1})  &= & \frac{\gamma b_k(t_{k-1} - t_k + 1)}{2t_k}, \vspace{1ex}\\
b_k\left(\frac{\eta_{k-1}t_k}{t_k-1} - \gamma_{k-1} \right) &= & 0.
\end{array}\right.
}
Using these expressions, we can simplify \eqref{eq:NAG_scheme10_proof7} as 
\myeqn{
\arraycolsep=0.2em
\begin{array}{lcl}
\Qc_k - \Qc_{k+1} \geq \frac{b_k(1 - L^2\gamma^2)t_k}{2L^2\gamma(t_k-1)} \norms{Gy_k - Gz_k}^2 + \frac{\gamma b_k(t_{k-1} - t_k + 1)}{2t_k}\norms{Gy_{k-1}}^2,
\end{array}
}
which proves the first estimate of \eqref{eq:NAG_scheme10_descent}.
}

\revise{
Finally, using the definition \eqref{eq:NAG_scheme10_Lyapunov} of $\Qc_k$, $x_k = y_{k-1} - \gamma_{k-1}Gy_{k-1}$, the monotonicity of $G$, and \eqref{eq:NAG_scheme10_param}, with a similar argument as in Lemma~\ref{le:NAG_scheme00_choice_of_params}, we can show that
\myeqn{
\arraycolsep=0.3em
\begin{array}{lcl}
\Qc_k  &\geq & \left( a_k - \frac{b_k^2}{4t_k^2} - \frac{b_k\gamma_{k-1}}{t_k}\right)\norms{Gy_{k-1}}^2 = \frac{\gamma b_k}{2t_k} \left( t_{k-1} - \frac{b_k}{2\gamma t_k} \right)\norms{Gy_{k-1}}^2,
\end{array} 
}
which proves  the second estimate of \eqref{eq:NAG_scheme10_descent}.
}
\Eproof
\end{proof}

Now, we can state the main convergence result of \eqref{eq:NAG_scheme10} in the following theorem.

\begin{theorem}\label{th:NAG_scheme10_convergence}
\revise{Suppose that $G$ in \eqref{eq:MI} is monotone and $L$-Lipschitz continuous.
Let $\sets{(x_k, y_k, z_k)}$ be generated by \eqref{eq:NAG_scheme10} using $t_k := k+\omega$ and \eqref{eq:NAG_scheme10_param} for a given $\omega > 1$.
Then, for all $k\geq 0$, we have 
\myeq{eq:NAG_scheme10_convergence1}{
\begin{array}{llcl}
& \norms{Gy_k}  & \leq &  \frac{2}{\gamma(k+\omega)} \norms{y_0 - y^{\star}}, \vspace{1ex}\\
& \sum_{l=0}^k(l+\omega)^2\norms{Gy_l-Gz_l}^2 & \leq & \frac{2L^2}{1-L^2\gamma^2}\norms{y_0 - y^{\star}}^2.
\end{array}
}
}
\end{theorem}

\begin{proof}
\revise{
Since $t_k := k + \omega$ for any $\omega > 1$, we have  $b_{k+1} = \frac{b_kt_{k+1}}{t_k-1} = \frac{b_k(k+\omega+1)}{k+\omega-1}$.
By induction, we get $b_k = \frac{b_1(k+\omega)(k+\omega-1)}{\omega(\omega+1)}$.
Using this expression and choosing $b_1 := \gamma\omega(\omega + 1)$, we obtain from the second line of \eqref{eq:NAG_scheme10_descent} the following bound:
\begin{equation}\label{eq:th5_proof100}
\Qc_{k+1} \geq \frac{\gamma b_1(k+\omega)^2}{2\omega(\omega+1)}\left(1 - \frac{b_1}{2\gamma\omega(\omega+1)}\right)\norms{Gy_k}^2 = \frac{\gamma^2(k+\omega)^2}{4} \norms{Gy_k}^2.
\end{equation}
Next, since $t_k := k + \omega$, we also have $t_{k-1} - t_k + 1 = 0$.
Moreover, since $0 < \gamma \leq \frac{1}{L}$ and $b_k = \gamma (k+\omega)(k+\omega-1)$, the first line of \eqref{eq:NAG_scheme10_descent} leads to 
\begin{equation}\label{eq:th5_proof101}
\Qc_k - \Qc_{k+1} \geq  \frac{(1-L^2\gamma^2)(k+\omega)^2}{2L^2}\norms{Gy_k - Gz_k}^2 \geq 0.
\end{equation}
However, since $x_0 = z_0$ and $\gamma Gy_{-1} = x_0 - y_{-1}$ from the first line of \eqref{eq:NAG_scheme10}, we have $\Qc_0 = \frac{a_0}{\gamma^2}\norms{x_0 - y_{-1}}^2 + \norms{x_0 - y^{\star} }^2$.
If we choose $y_{-1} := x_0 = y_0$, then we get $\Qc_0  = \norms{y_0 - y^{\star} }^2$.
From \eqref{eq:th5_proof101}, by induction, we obtain $\Qc_{k+1} \leq \Qc_0 = \norms{y_0 - y^{\star} }^2$.
Combining this inequality, and \eqref{eq:th5_proof100}, we get the first line of \eqref{eq:NAG_scheme10_convergence1}.
Finally, summing up \eqref{eq:th5_proof101} from $l := 0$ to $l := k$, we can deduce the second line of \eqref{eq:NAG_scheme10_convergence1}.
}
\Eproof
\end{proof}

\beforesubsec
\subsection{EAG for co-monotone case and its Nesterov's acceleration}
\aftersubsec
We consider the variant of EAG in \cite{lee2021fast} for the co-monotone operator $G$.
Recall that the operator $G$ in \eqref{eq:MI} is said to be $\rho$-comonotone if $\iprods{Gx - Gy, x - y} \geq \rho\norms{Gx - Gy}^2$ for all $x, y \in \R^p$, where $\rho < 0$.
In this subsection, we consider the case that $G$ is also $L$-Lipschitz continuous and $\rho$ satisfies the condition $ -\frac{1}{2L} < \rho \leq \frac{1}{L}$.
Hence, it covers three cases: co-coerciveness when $\rho > 0$, monotonicity when $\rho = 0$, and co-monotonicity when $\rho < 0$.
The second case has been studied in Subsection~\ref{subsec:EAG2NAG_scheme}.

More specifically, \cite{lee2021fast} proposes a variant of  \eqref{eq:EAG_scheme30} to solve \eqref{eq:MI} as follows:
\myeq{eq:p7:EAG_iter_com}{
\arraycolsep=0.3em
\left\{\begin{array}{lcl}
z_{k+1} &:= & \beta_ky_0 + (1-\beta_k)y_k - (1-\beta_k)\left(2\rho + \eta_k \right) Gy_k, \vspace{1ex}\\
y_{k+1} &:= & \beta_ky_0 + (1 - \beta_k)y_k - 2\rho(1 - \beta_k) Gy_k - \eta_kGz_{k+1},
\end{array}\right.
}
where $\beta_k := \frac{1}{k+1}$,  $\eta_k := \frac{1}{L}$.
If $\rho := 0$, then \eqref{eq:p7:EAG_iter_com} reduces to \eqref{eq:NAG_scheme10}.

As proved in \cite{lee2021fast}, the following convergence guarantee is obtained for \eqref{eq:p7:EAG_iter_com}:
\myeq{eq:p7:EAG_iter_com_convergence}{
\norms{Gy_k}^2 \leq \frac{4L^2\norms{y_0 - y^{\star}}^2}{(1 + 2\rho L)k^2}, \quad\forall k \geq 1.
}
Here, the key condition is $\rho > -\frac{1}{2L}$, which allows one to handle a class of  nonmonotone operators $G$.

Let us rewrite \eqref{eq:p7:EAG_iter_com} in a different form.
First, we have $y_{k+1} - z_{k+1} = -\eta_kGz_{k+1} + \eta_k(1-\beta_k)Gy_k$.
Hence, we get $2\rho(1-\beta_k)Gy_k = \frac{2\rho}{\eta_k}(y_{k+1} - z_{k+1}) + 2\rho Gz_{k+1}$.
Then, we have $y_{k+1} = \beta_ky_0 + (1-\beta_k)y_k - \eta_kGz_{k+1} -  \frac{2\rho}{\eta_k}(y_{k+1} - z_{k+1}) - 2\rho Gz_{k+1}$.
This expression implies that $(\eta_k + 2\rho)y_{k+1} = 2\rho z_{k+1} + \eta_k[ \beta_ky_0 + (1-\beta_k)y_k]  - \eta_k(\eta_k + 2\rho)Gz_{k+1}$.
Therefore, we eventually arrive at
\myeqn{
y_{k+1} = (1-\tau_k)z_{k+1} + \tau_k \big[  \beta_ky_0 + (1-\beta_k)y_k  - \tfrac{\eta_k}{\tau_k} Gz_{k+1} \big], \quad \text{where} \quad \tau_k := \tfrac{\eta_k}{\eta_k + 2\rho} > 0.
}
Overall, the scheme \eqref{eq:p7:EAG_iter_com} can be rewritten as
\myeq{eq:p7:EAG_iter_com_rewritten}{
\arraycolsep=0.3em
\left\{\begin{array}{lcl}
z_{k+1} &:= & \beta_ky_0 + (1-\beta_k)y_k - \frac{(1-\beta_k)\eta_k}{\tau_k} Gy_k, \vspace{1ex}\\
w_{k+1} &:= & \beta_ky_0 + (1 - \beta_k)y_k - \frac{\eta_k}{\tau_k}Gz_{k+1}, \vspace{1ex}\\
y_{k+1} &:= & (1-\tau_k)z_{k+1} + \tau_kw_{k+1}.
\end{array}\right.
}
The last line is a convex combination of the first two lines of the accelerated variant of EG+ from \cite{diakonikolas2021efficient}.
If $\tau_k = 1$ (i.e. $G$ is monotone), then \eqref{eq:p7:EAG_iter_com_rewritten} reduces to \eqref{eq:EAG_scheme30}.

Now, we derive Nesterov's accelerated interpretation of \eqref{eq:p7:EAG_iter_com}.
Following the same derivation as of \eqref{eq:NAG_scheme10}, we can show that \eqref{eq:p7:EAG_iter_com} is equivalent to
\myeq{eq:EAG_iter_com_Nest}{
\arraycolsep=0.3em
\left\{\begin{array}{lcl}
x_{k+1} & := & y_k -   (\eta_k + 2\rho) Gy_k, \vspace{1ex}\\
z_{k+1} & := & x_{k+1} +  \theta_k(x_{k+1} - x_k) + \nu_k (z_k  - x_{k+1}), \vspace{1ex}\\
y_{k+1} & := & z_{k+1} - \eta_k(Gz_{k+1} - (1-\beta_k)Gy_k). 
\end{array}\right.
}
This Nesterov's accelerated interpretation reduces to \eqref{eq:NAG_scheme10} when $\rho = 0$.
Moreover, it still achieves $\BigO{1/k^2}$ rate even when $G$ is co-monotone (i.e. $-\frac{1}{2L} < \rho < 0$), which is non-monotone.
The analysis in Theorem~\ref{th:NAG_scheme10_convergence} can be applied to \eqref{eq:EAG_iter_com_Nest}, but we omit it here.

\beforesubsec
\subsection{Nesterov's interpretation of the Halpern-type PEG method}\label{subsec:PEAG2NAG_scheme}
\aftersubsec
Our final step is to consider the Halpern-type past-extragradient (PEG) scheme for solving \eqref{eq:MI} studied in \cite{tran2021halpern}, which can be written as
\myeq{eq:PEAG_scheme00}{
\arraycolsep=0.3em
\left\{\begin{array}{lcl}
z_{k+1} &:= & \beta_ky_0 + (1-\beta_k)y_k - \eta_kGz_k, \vspace{1ex}\\
y_{k+1} &:= & \beta_ky_0 + (1-\beta_k)y_k - \hat{\eta}_kGz_{k+1}.
\end{array}\right.
}
Here, $z_0 := y_0$, $\beta_k \in (0, 1)$, and $\eta_k, \hat{\eta}_k > 0$ are given parameters, which will be determined in Theorem~\ref{th:NAG_scheme50_convergence}.
\revise{The convergence of \eqref{eq:PEAG_scheme00} has been proven in  \cite{tran2021halpern}.
However, we provide a slightly different variant in Theorem~\ref{th:NAG_scheme50_convergence} with $\eta_k \neq \hat{\eta}_k$ and with a simple parameter update for $\eta_k$ and $\hat{\eta}_k$ compared to  \cite{tran2021halpern}. 
We also provide a range of $\hat{\eta}_k$ instead of fixing it at $\hat{\eta}_k = \frac{1}{2L}$.
}

\begin{theorem}\label{th:NAG_scheme50_convergence}
\revise{Assume that $G$ in \eqref{eq:MI} is monotone and $L$-Lipschitz continuous and $\zer{G}\neq\emptyset$.
Let $\sets{(y_k, z_k)}$ be generated by \eqref{eq:PEAG_scheme00} to solve \eqref{eq:MI} using $\beta_k := \frac{1}{k+\omega}$, $\hat{\eta}_k := \hat{\eta} \in \left(0, \frac{1}{2L}\right]$, and $\eta_k := \hat{\eta}(1-\beta_k)$, where $\omega > 1$ is given.
Then, we have
\myeq{eq:NAG_scheme50_convergence}{
\arraycolsep=0.2em
\left\{ \begin{array}{lcl}
\norms{Gy_k}^2 + 2L^2\norms{z_k - y_k}^2 & \leq & \frac{C_0 \norms{y_0 - y^{\star}}^2}{(k+\omega - 1)^2}, \vspace{1ex}\\
\norms{Gz_k}^2 &  \leq &     \frac{3C_0 \norms{y_0 - y^{\star}}^2}{2(k+\omega - 1)^2}, \vspace{1ex}\\
\psi \cdot \sum_{k=1}^{\infty}(k+\omega-1)^2\big[\norms{Gy_k - Gy_{k-1}}^2 + L^2\norms{z_k - y_k}^2 \big] &\leq & C_0\norms{y_0 - y^{\star}}^2,
\end{array}\right.
}
where $C_0 := \frac{2\omega\hat{\eta}^2L^2[1 + (\omega-1)^2] + 4(\omega-1)}{\hat{\eta}^2(\omega-1)} > 0$ and $\psi := \frac{1-4L^2\hat{\eta}^2}{8L^2\hat{\eta}} \geq 0$.
}
\end{theorem}

\begin{proof}
\revise{We consider a Lyapunov function $\Lc_k := p_k\norms{Gy_k}^2 + q_k\iprods{Gy_k, y_k - y_0}$.
With a similar proof as in \cite{tran2021halpern}, by choosing $q_{k+1} = \frac{q_k}{1-\beta_k}$ and using \eqref{eq:PEAG_scheme00}, we have
\myeq{eq:PEAG_scheme00_proof2}{
\hspace{-4ex}
\arraycolsep=0.3em
\begin{array}{lcl}
\Lc_k - \Lc_{k+1} &\geq & \left( p_k - \frac{q_k\eta_k}{2\beta_k}\right)\norms{Gy_k}^2 + \left( \frac{q_k}{8L^2\eta_k\beta_k} - p_{k+1} \right)\norms{Gy_{k+1}}^2   \vspace{1ex}\\
&& + {~}  \frac{q_k(1 - 4L^2\hat{\eta}_k^2)}{8L^2\eta_k\beta_k}\norms{Gz_{k+1}}^2 + \frac{q_k}{\beta_k}\left( \frac{\hat{\eta}_k}{1-\beta_k}  - \frac{1}{4L^2\eta_k}\right) \iprods{Gy_{k+1}, Gz_{k+1}}  \vspace{1ex}\\
&& + {~}  \frac{q_k}{8\eta_k\beta_k}\norms{z_{k+1} - y_{k+1}}^2 - \frac{L^2q_k\eta_k}{2\beta_k}\norms{z_k - y_k}^2.
\end{array}
\hspace{-4ex}
}
Let us choose $\beta_k := \frac{1}{k+\omega}$ for some $\omega > 1$,  $\hat{\eta}_k := \hat{\eta} \in \left(0, \frac{1}{2L}\right]$, $\eta_k := \hat{\eta}_k(1-\beta_k) = \hat{\eta}(1-\beta_k)$, and $p_k := \frac{\hat{\eta}q_k}{2\beta_{k-1}}$.
Then, we have $q_{k+1} = \frac{q_k}{1-\beta_k} = \frac{q_k(k+\omega)}{k+\omega-1} = \frac{q_0(k+\omega)}{\omega}$.
Next,  using these choices of parameters and $q_{k+1} = \frac{q_k}{1-\beta_k}$, \eqref{eq:PEAG_scheme00_proof2} reduces to
\myeq{eq:PEAG_scheme00_proof2b}{ 
\arraycolsep=0.3em
\begin{array}{llcl}
\Ec_k &  \geq & \Ec_{k+1} + \frac{q_{k+1}(1 - 4L^2\hat{\eta}^2)(k+\omega)}{8L^2\hat{\eta}}\left[ \norms{Gy_{k+1} - Gy_k}^2  + L^2\norms{z_{k+1} - y_{k+1}}^2 \right],
\end{array}
}
where $\Ec_k   :=   \Lc_k +  \frac{\hat{\eta} L^2q_k(k+\omega-1)}{2}\norms{z_k - y_k}^2$.
}

\revise{Now, since $\Lc_k \geq \frac{p_k}{2}\norms{Gy_k}^2 - \frac{q_k^2}{2p_k}\norms{y_0 - y^{\star}}^2$ (see \cite{yoon2021accelerated}), we have
\myeq{eq:PEAG_scheme00_proof2c}{ 
\arraycolsep=0.3em
\begin{array}{lcl}
\Ec_k &\geq &  \frac{\hat{\eta}q_k(k+\omega-1)}{4}\norms{Gy_k}^2 + \frac{\hat{\eta} L^2q_k(k+\omega-1)}{2}\norms{z_k - y_k}^2 - \frac{q_k}{\hat{\eta}(k+\omega-1)}\norms{y_0 - y^{\star}}^2 \vspace{1ex}\\
& = & \frac{\hat{\eta}q_0(k+\omega-1)^2}{4\omega}\big[ \norms{Gy_k}^2 + 2L^2\norms{y_k - z_k}^2 \big] - \frac{q_0}{\hat{\eta}\omega}\norms{y_0 - y^{\star}}^2.
\end{array}
}
If $2L\hat{\eta} \leq 1$, then from \eqref{eq:PEAG_scheme00_proof2b}, we have $\Ec_{k+1} \leq \Ec_k$.
By induction, it leads to $\Ec_k \leq \Ec_0$.
Furthermore, by the Lipschitz continuity of $G$ and $Gy^{\star} = 0$, we have $\norms{Gy_0} \leq L\norms{y_0 - y^{\star}}$.
Using this estimate, we can easily show that
\myeqn{ 
\arraycolsep=0.3em
\begin{array}{lcl}
\Ec_0 & = & \frac{\hat{\eta} q_0}{2(\omega - 1)}\norms{Gy_0}^2  + \frac{\hat{\eta}L^2 q_0(\omega-1)}{2}\norms{y_0 - y^{\star}}^2 \leq \frac{\hat{\eta}L^2 q_0[1 + (\omega-1)^2]}{2(\omega-1)}\norms{y_0 - y^{\star}}^2.
\end{array}
}
Hence, combining this bound and  $\Ec_k \leq \Ec_0$, we get $\Ec_k \leq  \frac{\hat{\eta}L^2 q_0[1 + (\omega-1)^2]}{2(\omega-1)}\norms{y_0 - y^{\star}}^2$.
Utilizing \eqref{eq:PEAG_scheme00_proof2c}, the last inequality leads to the first line of \eqref{eq:NAG_scheme50_convergence}.
}

\revise{By Young's inequality and the Lipschitz continuity of $G$, we have  $\norms{Gz_k}^2 \leq \frac{3}{2}\norms{Gy_k}^2 + 3\norms{Gz_k - Gy_k}^2  \leq \frac{3}{2} \left[ \norms{Gy_k}^2 + 2 L^2\norms{z_k - y_k}^2  \right]$.
Combining this inequality and the first line of \eqref{eq:NAG_scheme50_convergence}, we get the second line of \eqref{eq:NAG_scheme50_convergence}.
}

\revise{Finally, summing up \eqref{eq:PEAG_scheme00_proof2c} from $k :=0$ to $k := K$, and then using \eqref{eq:PEAG_scheme00_proof2c} and the upper bound of $\Ec_0$, we have
\myeqn{
\arraycolsep=0.3em
\begin{array}{ll}
\frac{q_0(1 - 4L^2\hat{\eta}^2)}{8L^2\hat{\eta}} & \sum_{k=0}^K(k+\omega)^2 \left[ \norms{Gy_{k+1} - Gy_k}^2  + L^2\norms{z_{k+1} - y_{k+1}}^2 \right] \leq \Ec_0 - \Ec_{K+1} \vspace{1ex}\\
&\leq \frac{\hat{\eta}L^2 q_0[1 + (\omega-1)^2]}{2(\omega-1)}\norms{y_0 - y^{\star}}^2 + \frac{q_0}{\hat{\eta}\omega}\norms{y_0 - y^{\star}}^2.
\end{array}
}
Simplifying this inequality and letting $K\to\infty$, we obtain the third line of \eqref{eq:NAG_scheme50_convergence}.
\Eproof
}
\end{proof}

\revise{
Note that if we choose $\omega := 2$  in Theorem~\ref{th:NAG_scheme50_convergence}, then $\beta_k = \frac{1}{k+2}$ and $C_0 =  \frac{4(2L^2\hat{\eta}^2 + 1)}{\hat{\eta}^2}$.
This constant factor is larger than the one in \eqref{eq:NAG_scheme10_convergence1} of Theorem~\ref{th:NAG_scheme10_convergence}.
However, \eqref{eq:PEAG_scheme00} only requires one evaluation of $G$ per iteration compared to two evaluations as in  \eqref{eq:EAG_scheme30}.
If $0 < \hat{\eta} < \frac{1}{2L}$, then the last  summable bound in the third line of \eqref{eq:NAG_scheme50_convergence} is not vanished.}

By following the same arguments as \eqref{eq:NAG_scheme10}, we can derive Nesterov's accelerated interpretation of \eqref{eq:PEAG_scheme00}.
This result is stated in the following theorem.

\begin{theorem}\label{th:HP2Nes_scheme_PEAG}
Let $\beta_k \in (0, 1)$ and $\eta_k, \hat{\eta}_k > 0$ be given as in \eqref{eq:PEAG_scheme00}.
Let $\sets{(\hat{x}_k, z_k)}$ be generated by the following scheme:
\myeq{eq:NAG_scheme50}{
\arraycolsep=0.2em
\left\{\begin{array}{lcl}
\hat{x}_{k+1} &:= & z_k - \hat{\gamma}_kGz_k, \vspace{1ex}\\
z_{k+1} & := & \hat{x}_{k+1}  +  \theta_k(\hat{x}_{k+1} - \hat{x}_k) +  \nu_k(z_k - \hat{x}_{k+1}) + \kappa_k(z_{k-1} - \hat{x}_k)  \vspace{1ex}\\
&& - {~} \zeta_k(z_{k-2} - \hat{x}_{k-1}),
\end{array}\right.
}
starting from $z_{-2} = z_{-1} = z_0 = \hat{x}_{-1} = \hat{x}_0 := y_0$, and using $\theta_k := \frac{\beta_k(1-\beta_{k-1})}{\beta_{k-1}}$, $\nu_k := \frac{\beta_k}{\beta_{k-1}}$, $\hat{\gamma}_k :=  \hat{\eta}_{k-1} +  \frac{\eta_k}{1-\beta_k}$, $\kappa_k :=  \frac{\eta_{k-1}(1 - \beta_k)}{\hat{\gamma}_{k-1}} $ and $\zeta_k := \frac{\theta_k\eta_{k-2}}{\hat{\gamma}_{k-2}}$ provided that $\hat{\eta}_{-2} = \hat{\eta}_{-1} = \hat{\eta}_0 =  \eta_{-2} = \eta_{-1} =  \eta_0$.
Let $y_k :=  z_k - \hat{\eta}_{k-1}Gz_k + \eta_{k-1}Gz_{k-1}$.
Then, $\sets{(y_k, z_k)}$  is identical to the one generated by \eqref{eq:PEAG_scheme00} starting from $y_0$.
\end{theorem}

\begin{proof}
Similar to \eqref{eq:NAG_scheme10}, the following scheme is equivalent to \eqref{eq:PEAG_scheme00}:
\myeq{eq:NAG_scheme30}{
\arraycolsep=0.3em
\left\{\begin{array}{lcl}
x_{k+1} &:= & y_k - \gamma_k Gz_k, \vspace{1ex}\\
z_{k+1} & := & x_{k+1}  +  \theta_k(x_{k+1} - x_k) +  \nu_k(z_k - x_{k+1}), \vspace{1ex}\\
y_{k+1} & := & z_{k+1} - \hat{\eta}_kGz_{k+1} + \eta_kGz_k,
\end{array}\right.
}
where $\gamma_k := \frac{\eta_k}{1-\beta_k}$, $\theta_k := \frac{\beta_k(1-\beta_{k-1})}{\beta_{k-1}}$, and $\nu_k := \frac{\beta_k}{\beta_{k-1}}$.

Next, from the third line of \eqref{eq:NAG_scheme30}, we have $y_k = z_k - \hat{\eta}_{k-1}Gz_k + \eta_{k-1}Gz_{k-1}$.
Hence, we can eliminate $y_k$ in \eqref{eq:NAG_scheme30}  to get $x_{k+1} = z_k - (\hat{\eta}_{k-1} + \gamma_k)Gz_k + \eta_{k-1}Gz_{k-1}$.
In this case, \eqref{eq:NAG_scheme30} can be written equivalently to 
\myeq{eq:NAG_scheme50_inter}{
\arraycolsep=0.3em
\left\{\begin{array}{lcl}
x_{k+1} &:= & z_k - (\hat{\eta}_{k-1} + \gamma_k)Gz_k + \eta_{k-1}Gz_{k-1}, \vspace{1ex}\\
z_{k+1} & := & x_{k+1}  +  \theta_k(x_{k+1} - x_k) +  \nu_k(z_k - x_{k+1}).
\end{array}\right.
}
This scheme can be viewed as Nesterov's accelerated interpretation of \eqref{eq:PEAG_scheme00}.
However, if we denote $\hat{x}_{k+1} := z_k - (\hat{\eta}_{k-1} + \gamma_k)Gz_k$, then $x_{k+1} = \hat{x}_{k+1} +  \eta_{k-1}Gz_{k-1}$.
Substituting this expression into the second line of \eqref{eq:NAG_scheme50_inter}, we get 
\myeqn{
z_{k+1} = \hat{x}_{k+1}  +  \theta_k(\hat{x}_{k+1} - \hat{x}_k) +  \nu_k(z_k - \hat{x}_{k+1}) + \eta_{k-1}(1+\theta_k - \nu_k)Gz_{k-1} - \theta_k\eta_{k-2}Gz_{k-2}.
}
Since $Gz_{k-1} = \frac{1}{\hat{\eta}_{k-2} + \gamma_{k-1}}(z_{k-1} - \hat{x}_k)$, we can further write \eqref{eq:NAG_scheme50_inter} as
\myeqn{ 
\arraycolsep=0.2em
\left\{\begin{array}{lcl}
\hat{x}_{k+1} &:= & z_k - \hat{\gamma}_kGz_k, \vspace{1ex}\\
z_{k+1} & := & \hat{x}_{k+1}  +  \theta_k(\hat{x}_{k+1} - \hat{x}_k) +  \nu_k(z_k - \hat{x}_{k+1}) \vspace{1ex}\\
&& + {~} \kappa_k(z_{k-1} - \hat{x}_k) - \zeta_k(z_{k-2} - \hat{x}_{k-1}),
\end{array}\right.
}
where $\hat{\gamma}_k :=  \hat{\eta}_{k-1} + \gamma_k$, $\kappa_k := \frac{\eta_{k-1}(1+\theta_k - \nu_k)}{\hat{\eta}_{k-2} + \gamma_{k-1}} =  \frac{\eta_{k-1}(1 - \beta_k)}{\hat{\gamma}_{k-1}} $ and $\zeta_k :=  \frac{\theta_k\eta_{k-2}}{\hat{\eta}_{k-3} + \gamma_{k-2}} = \frac{\theta_k\eta_{k-2}}{\hat{\gamma}_{k-2}}$.
The last scheme is exactly \eqref{eq:NAG_scheme50}.
\Eproof
\end{proof}

Clearly, the new Nesterov's accelerated  scheme \eqref{eq:NAG_scheme50} for solving \eqref{eq:MI} has three correction terms instead of two as in \eqref{eq:NAG_scheme00_2corr}.
Since our transformation is equivalent, the convergence of \eqref{eq:NAG_scheme50}  is still guaranteed by Theorem~\ref{th:NAG_scheme50_convergence}.
More specifically, if we choose $\beta_k$, $\hat{\eta}_k$, and $\eta_k$ as in Theorem~\ref{th:NAG_scheme50_convergence}, then we obtain the following corollary.

\begin{corollary}\label{co:eq:NAG_scheme50_convergence}
\revise{Assume that $G$ in \eqref{eq:MI} is monotone and $L$-Lipschitz continuous and $\zer{G}\neq\emptyset$.
Let $\sets{(\hat{x}_k, z_k)}$ be generated by \eqref{eq:NAG_scheme50} to solve \eqref{eq:MI} using $\hat{\gamma}_k := \hat{\gamma} \in \left(0,  \frac{1}{L}\right]$ for all $k\geq 0$, and
\myeqn{ 
\arraycolsep=0.3em
\begin{array}{ll}
\theta_k := \frac{k+\omega-2}{k+\omega}, \quad \nu_k := \frac{k+\omega-1}{k+\omega},  \quad \kappa_k := \frac{k + \omega - 2}{2(k+\omega)}, \quad \text{and} \quad \zeta_k := \begin{cases} 0 &\text{if $k=0$,} \\  \frac{k + \omega -3}{2(k+\omega)} &\text{if $k\geq 1$},\end{cases}
\end{array}
}
where $\omega > 1$ is a given parameter.
Then, we have
\myeq{eq:NAG_scheme50_convergence2}{
\arraycolsep=0.3em
\left\{ \begin{array}{lcl}
\norms{G\hat{x}_{k+1}}^2  & \leq &  \frac{3(1 + L^2\hat{\gamma}^2)\tilde{C}_0 \norms{y_0 - y^{\star}}^2}{(k+\omega - 1)^2}, \vspace{1ex}\\
\norms{Gz_k}^2 &  \leq &   \frac{3\tilde{C}_0 \norms{y_0 - y^{\star}}^2}{2(k+\omega - 1)^2},
\end{array}\right.
}
where $\tilde{C}_0 := \frac{2\omega\hat{\gamma}^2L^2[1 + (\omega-1)^2] + 16(\omega-1)}{\hat{\gamma}^2(\omega-1)} > 0$ is rendered from $C_0$ of  Theorem~\ref{th:NAG_scheme50_convergence}.
}
\end{corollary}

\begin{proof}
\revise{First, by the choice of $\beta_k := \frac{1}{k+\omega}$ for some $\omega > 1$, $\hat{\eta}_k := \hat{\eta} \in \left(0,  \frac{1}{2L}\right]$, and $\eta_k := \hat{\eta}(1-\beta_k)$ in Theorem \ref{th:NAG_scheme50_convergence}, 
using the update rules of parameters in Theorem~\ref{th:HP2Nes_scheme_PEAG}, we can easily  show that $\theta_k =  \frac{k+\omega-2}{k+\omega}$, $\nu_k = \frac{k+\omega-1}{k+\omega}$, $\hat{\gamma}_k = \hat{\gamma} := 2\hat{\eta}$, $\kappa_k =  \frac{k+\omega-2}{2(k+\omega)}$, and $\zeta_k = \frac{k+\omega-3}{2(k+\omega)}$ if $k\geq 1$, and $\zeta_k = 0$ if $k=0$, as given in Corollary~\ref{co:eq:NAG_scheme50_convergence}.
}

\revise{Next, since  $\hat{x}^{k+1} = z_k - \hat{\gamma}_kGz_k = z_k - \hat{\gamma} Gz_k$ due to \eqref{eq:NAG_scheme50}, we have $\norms{G\hat{x}_{k+1}}^2 \leq 2\norms{G\hat{x}_{k+1} - Gz_k}^2 + 2\norms{Gz_k}^2 \leq 2L^2\norms{\hat{x}_{k+1} - z_k}^2 + 2\norms{Gz_k}^2 = 2L^2\hat{\gamma}^2\norms{Gz_k}^2 + 2\norms{Gz_k}^2 = 2(1 + L^2\hat{\gamma}^2)\norms{Gz_k}^2$.
Combining this inequality and the second line of \eqref{eq:NAG_scheme50_convergence}, and noting that $\hat{x}_0 = y_0$, we obtain the first line of \eqref{eq:NAG_scheme50_convergence2}.
Finally, the second line of  \eqref{eq:NAG_scheme50_convergence2} directly comes from the second line of  \eqref{eq:NAG_scheme50_convergence}.
}
\Eproof
\end{proof}

Unlike \eqref{eq:NAG_scheme00_2corr}, the new scheme \eqref{eq:NAG_scheme50} has convergence without the co-coerciveness of $G$.
It only requires $G$ to be monotone and $L$-Lipschitz continuous, and one evaluation of $G$ per iteration.

\beforesec
\section{Numerical Experiments}\label{sec:num_experiments}
\aftersec
In this section, we illustrate our theoretical results through two numerical examples.
The first one is to test Nesterov's accelerated variant \eqref{eq:NAG_scheme00} with two sets of parameters stated in Corollary~\ref{co:HP_scheme00_convergence} and Theorem~\ref{th:NAG_scheme00_convergence2}, respectively.
The second example is to test two Nesterov's accelerated schemes \eqref{eq:NAG_scheme10} and \eqref{eq:NAG_scheme50} of EAG and PEAG, respectively in Section~\ref{sec:EAG_sec}.

\beforesubsec
\subsection{The performance of Nesterov's accelerated variants of \eqref{eq:HP_scheme00}}\label{subsec:exam1}
\aftersubsec
We consider a linear regression model $b = Py + \epsilon$, and its corresponding least-squares problem, where $P \in \R^{n\times p}$, $y$ is a model parameter, and $\epsilon$ is a Gaussian noise of zero mean and variance $\sigma^2$. 
We define the operator $G$ from the normal equation of this least-squares problem, which is written as $Gy = P^{\top}(Py - b)$.
It is obvious to show that $Gy$ is $\frac{1}{L}$-co-coercive with $L = \norms{P^{\top}P}$.

We implement two variants of \eqref{eq:NAG_scheme00} with two sets of parameters stated in Corollary~\ref{co:HP_scheme00_convergence} and Theorem~\ref{th:NAG_scheme00_convergence2}, respectively to solve $Gy^{\star} = 0$ (i.e. $P^{\top}(Py^{\star} - b) = 0$).
The input data is generated as follows.
We choose two cases: Case 1 with $(n, p) = (500, 1000)$ and Case 2 with $(n, p) = (1000, 1000)$.
We generate $P$ randomly using the standard Gaussian distribution and then normalize  it to get unit columns.
Next, we generate $b := Py^{\natural} + \Nc(0, 0.1)$, where $y^{\natural}$ is a given vector generated from the standard Gaussian distribution, and $\Nc(0, 0.1)$ is a Gaussian noise of zero mean and variance $\sigma^2 = 0.1$.
This procedure makes sure that $y^{\star}$ exists.

\begin{figure}[ht!] 
\vspace{-3ex}
\begin{center}
\includegraphics[width=0.495\textwidth]{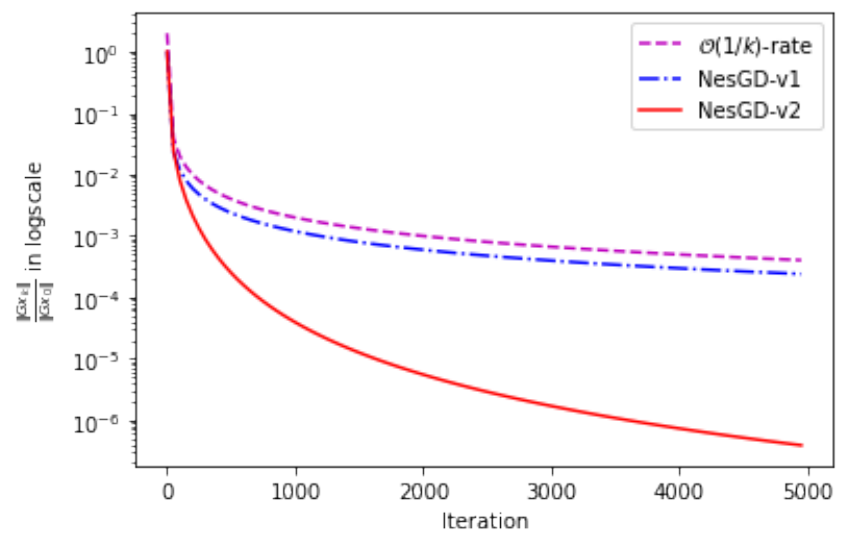}
\includegraphics[width=0.495\textwidth]{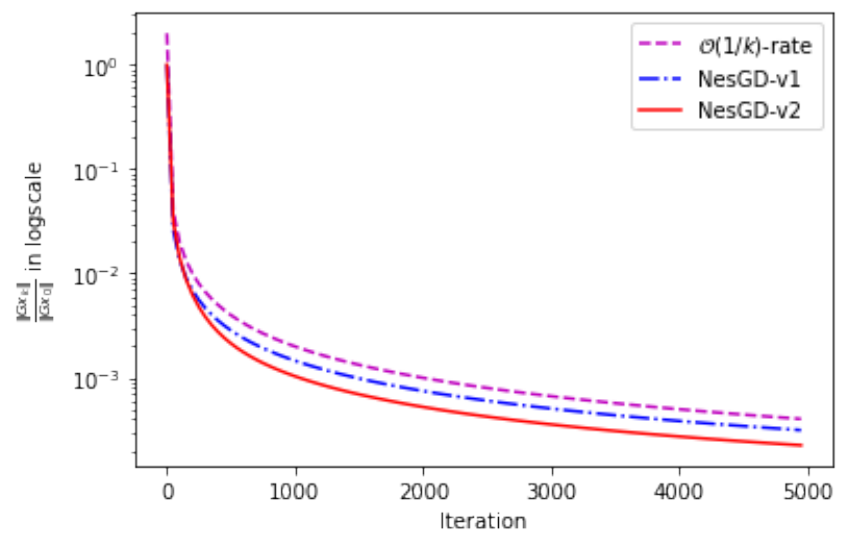}
\vspace{-3ex}
\caption{The convergence behavior of the two Nesterov's accelerated variants on two problem instances. Left: Case 1 with $(n, p) = (500, 1000)$, and Right: Case 2 with $(n, p) = (1000, 1000)$.}\label{fig:exam1}
\vspace{-4ex}
\end{center}
\end{figure}

We run both variants: \texttt{NesGD-v1} (Corollary~\ref{co:HP_scheme00_convergence}) and \texttt{NesGD-v2} (Theorem~\ref{th:NAG_scheme00_convergence2}) on these two instances up to $5000$ iterations.
\revise{For \texttt{NesGD-v2}, we choose $\omega := 3$ to update its parameters.}
The results of both variants are revealed in Figure~\ref{fig:exam1}, where the $y$-axis shows the relative error $\frac{\norms{Gx_k}}{\norms{Gx_0}}$ in log-scale.

Figure~\ref{fig:exam1} shows that \texttt{NesGD-v1} follows $\BigO{1/k}$-convergence rate in both cases, while \texttt{NesGD-v2} highly outperforms its $\BigO{1/k}$ theoretical rate  in Case 1.
In Case 2, \texttt{NesGD-v2} still follows its $\BigO{1/k}$-convergence rate.
Moreover, $\norms{Gx_k}$ is decreasing smoothly without oscillation w.r.t. the iteration counter $k$ in both variants.

\beforesubsec
\subsection{The performance of Nesterov's variants of EAG and PEAG}\label{subsec:exam2}
\aftersubsec
Our next experiment is to test Nesterov's accelerated variants \eqref{eq:NAG_scheme10} of EAG and \eqref{eq:NAG_scheme50} of PEAG.
For this purpose, we consider a minimax problem and the corresponding saddle operator $G$ as
\myeq{eq:exam_func1}{
\min_{u\in\R^m}\max_{v \in \R^n}\Big\{ \Hc(u, v) := f(u) + \iprods{Ku, v} - g(v) \Big\}, \ \text{and} \ 
Gy := \begin{bmatrix}\nabla{f}(u) + K^{\top}v\\ \nabla{g}(v) - Ku\end{bmatrix},
}
where $K \in \R^{m\times n}$, and $\nabla{f}$ and $\nabla{g}$ are the gradients of $L$-smooth and convex functions $f$ and $g$, respectively.
In this case, it is well-known that $G$ is monotone and Lipschitz continuous with $L := \sqrt{2}\big[\max\sets{L_f^2, L_g^2} + \norms{K}^2\big]^{1/2}$, where $L_f$ and $L_g$ are the Lipschitz constants of $\nabla{f}$ and $\nabla{g}$, respectively.

In our experiments, we choose $f(u) := \lambda\sum_{i=1}^{m}\ell_{\epsilon}(u_i)$, and $g(v) := \rho\sum_{i=1}^{n}\ell_{\epsilon}(v_i)$, where $\ell_{\epsilon}(\tau)$ is the well-known Huber function (i.e. $\ell_{\epsilon}(\tau) = \epsilon\vert\tau\vert - \frac{\epsilon^2}{2}$ if $\vert \tau \vert \geq \epsilon$, and $\ell_{\epsilon}(\tau) = \frac{\tau^2}{2}$, otherwise) and $\lambda > 0$ and $\rho > 0$ are given parameters. 
Clearly, one can easily show that $\nabla{f}(u) = \lambda(\ell'_{\epsilon}(u_1), \cdots, \ell'_{\epsilon}(u_m))^{\top}$ and $\nabla{g}(v) = \rho(\ell'_{\epsilon}(v_1), \cdots, \ell'_{\epsilon}(v_n))^{\top}$, where $\ell'_{\epsilon}(\tau) = \tau$ if $\vert \tau \vert < \epsilon$, $\ell'_{\epsilon}(\tau) = -\epsilon$ if $\tau \leq -\epsilon$, and $\ell'_{\epsilon}(\tau) = \epsilon$ if $\tau \geq \epsilon$.
To obtain data for our experiments, we generate a random matrix $K$ in $\R^{m\times n}$ from the standard Gaussian distribution and normalize its columns.
Then, we compute $\norms{K}$ and then choose $\lambda = \rho := \norms{K}$ and $\epsilon := 0.05$.

\begin{figure}[ht!] 
\vspace{-3ex}
\begin{center}
\includegraphics[width=0.495\textwidth]{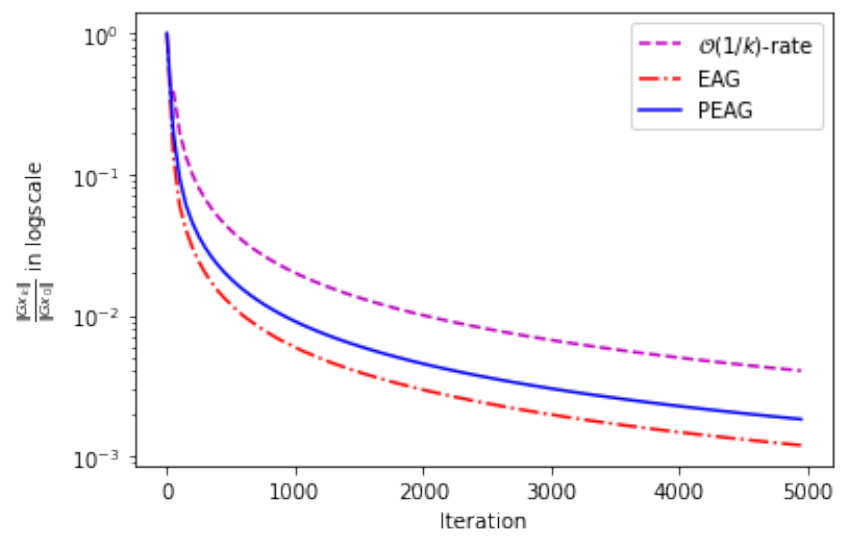}
\includegraphics[width=0.495\textwidth]{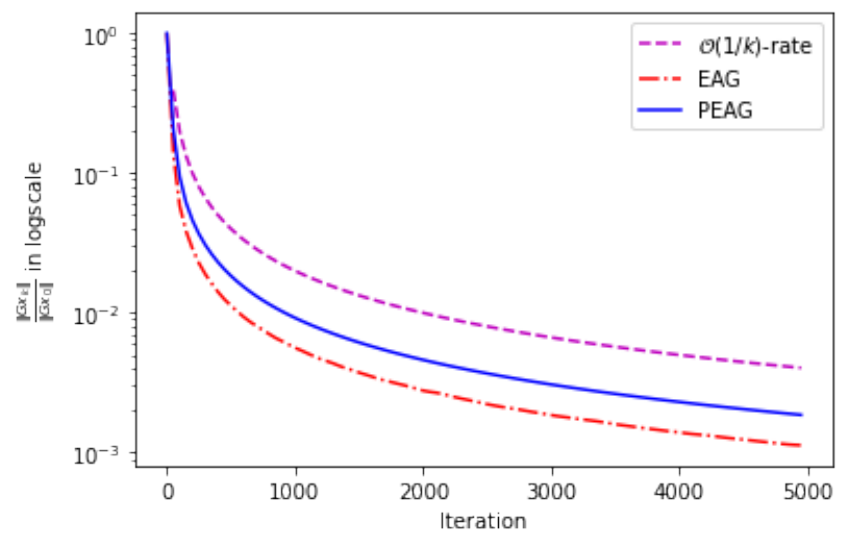}
\vspace{-3ex}
\caption{The performance of Nesterov's accelerated variants of EAG and PEAG on \eqref{eq:exam_func1}. Left: Case 1, and Right: Case 2}\label{fig:exam2}
\end{center}
\vspace{-4ex}
\end{figure}

We implement both Nesterov's accelerated variants \eqref{eq:NAG_scheme10} and \eqref{eq:NAG_scheme50} to solve problem \eqref{eq:exam_func1}.
The parameters are updated exactly as in Theorems~\ref{th:NAG_scheme10_convergence} and Corollary~\ref{co:eq:NAG_scheme50_convergence}, respectively \revise{with $\omega := 2$}. 
We run these two algorithms on two cases: Case 1: $(m, n) = (1000, 750)$ and Case 2: $(m, n) = (1000, 1000)$.
We set the number of iterations $k_{\max}$ at $k_{\max} := 5000$.
The convergence behavior of both schemes on $\frac{\norms{Gx_k}}{\norms{Gx_0}}$ is plotted in Figure~\ref{fig:exam2} for the two instances.

Figure~\ref{fig:exam2} shows that both schemes match very well the $\BigO{1/k}$-rate in terms of $\norms{Gx_k}$ as proven by our theory. 
Here, \eqref{eq:NAG_scheme10}  performs slightly better than \eqref{eq:NAG_scheme50}  due to its bigger stepsize. 
However, \eqref{eq:NAG_scheme50}   only requires one evaluation of $G$ at each iteration as opposed to two as in \eqref{eq:NAG_scheme10}.

{
\vspace{2ex}
\noindent\textbf{Data availability.}
The author confirms that all data used in this paper is generated synthetically.
The method for generating data is also described in the paper.
}

\vspace{2ex}
\noindent\textbf{Acknowledgements.}
This paper is based upon work partially supported by the National Science Foundation (NSF), grant no. NSF-RTG DMS-2134107 and the Office of Naval Research (ONR), grant No. N00014-20-1-2088 (2020-2023) and grant No. N00014-23-1-2588 (2023-2026).
We would like to express our sincere gratitude to the two anonymous reviewers for their invaluable comments and feedback to improve the paper.

\vspace{2ex}
\noindent\textbf{Conflicts of interest.}
The author declares that he has no conflict of interest of any kind.

%

\end{document}